\documentclass[reqno]{amsart}

\usepackage{lineno,hyperref}

\usepackage{graphicx,epstopdf,epsfig}
\usepackage{amsfonts,epsfig,fancyhdr,graphics, hyperref,amsmath,amssymb, amsthm}
\usepackage[inline]{enumitem}
\usepackage{mathtools}
\usepackage{tikz}
\usepackage{xcolor}
\usepackage{cleveref}
\usepackage{bm}
\usepackage{footmisc}

\makeatletter
\newcommand*\rel@kern[1]{\kern#1\dimexpr\macc@kerna}
\newcommand*\widebar[1]{%
	\begingroup
	\def\mathaccent##1##2{%
		\rel@kern{0.8}%
		\overline{\rel@kern{-0.8}\macc@nucleus\rel@kern{0.2}}%
		\rel@kern{-0.2}%
	}%
	\macc@depth\@ne
	\let\math@bgroup\@empty \let\math@egroup\macc@set@skewchar
	\mathsurround\z@ \frozen@everymath{\mathgroup\macc@group\relax}%
	\macc@set@skewchar\relax
	\let\mathaccentV\macc@nested@a
	\macc@nested@a\relax111{#1}%
	\endgroup
}
\makeatother

\makeatletter
\newcommand{\inlineitem}[1][]{%
	\ifnum\enit@type=\tw@
	{\descriptionlabel{#1}}
	\hspace{\labelsep}%
	\else
	\ifnum\enit@type=\z@
	\refstepcounter{\@listctr}\fi
	\quad\@itemlabel\hspace{\labelsep}%
	\fi}
\makeatother

\makeatletter
\g@addto@macro{\endabstract}{\@setabstract}
\newcommand{\authorfootnotes}{\renewcommand\thefootnote{\@fnsymbol\c@footnote}}%
\makeatother

\theoremstyle{plain}
\newtheorem{theorem}{Theorem}[section]
\newtheorem{corollary}[theorem]{Corollary}
\newtheorem{proposition}[theorem]{Proposition}
\newtheorem{lemma}[theorem]{Lemma}

\theoremstyle{definition}
\newtheorem{definition}[theorem]{Definition}
\newtheorem{remark}[theorem]{Remark}
\newtheorem{example}[theorem]{Example}
\newtheorem{observation}[theorem]{Observation}
\newtheorem{question}[theorem]{Question}

\newcommand{\be}{\mathbf e}
\newcommand{\cB}{\mathcal B}

\newcommand{\cT}{\mathcal T}

\usepackage[export]{adjustbox}
\usepackage{forest}
\usetikzlibrary{shapes.geometric,decorations.markings}
\forestset{
	my tree/.style={
		before typesetting nodes={
			for tree={
				if={>O+tt=!O_=&{content}{}{n}{1}}{
					label/.process={Ow}{content}{left:##1},
				}{
					label/.process={Ow}{content}{right:##1},
				},
				content=,
			},
		},
		where level=0{
			draw,
			baseline,
			fill,
			circle,
		}{
			draw,
			fill,
			circle,
		},
		for tree={inner sep=1.3pt, l sep+=1pt, s sep'+=15pt},
	},
	default preamble={my tree},
	cut/.style={
		tikz+={
			\path ()  -- (!u);
		}
	}
}

\begin{document}
	
	\begin{center}
		\LARGE 
		Kemeny's constant for a graph with bridges \par \bigskip
		
		\normalsize
		\authorfootnotes
		Jane Breen\footnote{Faculty of Science, Ontario Tech University, Oshawa, Ontario, Canada}, Emanuele Crisostomi\footref{foot1}, Sooyeong Kim\footnote{Department of Energy, Systems, Territory and Constructions Engineering, University of Pisa, 56121, Pisa, Italy.\label{foot1} (Contact: kswim2502@gmail.com)} \par \bigskip
		
		\today
	\end{center}
\begin{abstract}
	In this paper, we determine a formula for Kemeny's constant for a graph with multiple bridges, in terms of quantities that are inherent to the subgraphs obtained upon removal of all bridges and that can be computed independently. With the formula, we consider several optimization problems for Kemeny's constant for graphs with bridges, and we remark on the computational benefit of this formula for the computation of Kemeny's constant. Finally, we discuss some potential applications.
\end{abstract}



\section{Introduction}

\subsection{Motivation}
Kemeny's constant is an interesting quantifier of a finite discrete-time Markov chain, measuring the expected time it takes to travel between randomly-chosen states of the chain, where each state is selected with probability proportional to the long-run prevalence of the state. Originally defined in \cite{KemenySnell} as the expected time to reach a randomly-chosen state from a fixed starting state, this quantity is astonishingly independent of the choice of initial state, and thus is named Kemeny's constant. Since its introduction in \cite{KemenySnell}, this quantity has fascinated a number of researchers interested in providing an intuitive explanation of its constant property (see, for example, \cite{bini2018kemeny, Doyle2009, kirkland2021directed, Levene:KemenyInterpre}). Due to its interpretation in terms of how `well-mixing' a Markov chain is, it is not surprising that Kemeny's constant has been recently used for a wide range of applications. Some significant examples include the modeling of road traffic networks \cite{Crisostomi:Google, Salman2018}; city-scale models of pollution from tyre particles \cite{Singer2021}; vaccination and testing strategies in epidemiological networks \cite{Emma:Kemeny}; and automated stochastic surveillance strategies \cite{Bullo2015}.

Of particular interest is the study of the concept of Kemeny's constant for simple random walks on graphs, which are a subclass of Markov chains in which the vertices of the graph $G$ form the state space, and one models stochastic movement between adjacent vertices in the graph. As such, Kemeny's constant for such a Markov chain (denoted $\kappa(G)$) provides a global quantifier of how `well-connected' the vertices of the graph are, in terms of how long it takes to get from a randomly-chosen initial vertex to a randomly-chosen destination vertex. Large values of $\kappa(G)$ can indicate clustering in the graph, while small values of $\kappa(G)$ indicate that the graph $G$ has good expansion properties.
In addition, Kemeny's constant is closely connected to a number of many other well-known graph quantities which are commonly used in different scientific communities, such as the effective graph resistance \cite{wang2017kemeny}, or the subdominant eigenvalues \cite{Dudkina2021}.

Despite the many applications, there are still many open problems involving Kemeny's constant. In particular, one of the main challenges is determining how to compute it efficiently in large-scale networks, and how to update its value under minor modifications to the topology of the underlying graph (e.g., vertex or edge removal, vertex or edge addition in the graph \cite{Shimbo2016}), without having to recompute it from scratch.
This is exactly the objective of this manuscript, as we investigate the problem of how to compute Kemeny's constant for a random walk in a graph using a divide-and-conquer strategy, given the knowledge of Kemeny's constants for smaller subgraphs of the graph. We focus on the case that the graph $G$ has bridges (i.e.~edges whose deletion disconnects the graph into smaller components).

While this scenario is a special case, we believe that the results pave the way for future research along the same lines, leading to an easier recalculation of Kemeny's constant, and also for parallel computations of its value in some specific cases. 
Some work has already begun in this area. In \cite{altafini2022centrality}, the authors determine an expression for the change in Kemeny's constant after edge deletion in a connected graph. In \cite{faught20221}, a formula for Kemeny's constant is provided for graphs with cut vertices. In addition to providing a new expression for Kemeny's constant for graphs with bridges, we also improve on this previous work by exploring the technical expressions for computing Kemeny's constant and expressing them in terms of meaningful Markov chain quantities (such as accessibility indices and mean first passage times) in order to develop further intuition on the underlying meaning of Kemeny's constant. 

\subsection{Preliminaries}

A graph $G$ consists of a set of vertices $V(G)$ and a set of edges $E(G) \subseteq \{\{u, v\} \mid u, v \in V(G)\}$. We use $m_G$ to denote the number of edges. Two vertices $u$ and $v$ are said to be \emph{adjacent} if there is an edge $\{u, v\}$ in $E(G)$; this is also denoted by $u\sim v$. The \emph{degree} of a vertex $v$, denoted $\deg_G(v)$, is the number of vertices adjacent to $v$ in $G$. For $v,w\in V(G)$, we use $\mathrm{dist}_G(v,w)$ to denote the distance between $v$ and $w$; that is, the length of the shortest path between $v$ and $w$. A graph is said to be \emph{connected} if there is a path from $u$ to $v$ for every pair of vertices $u,v$ in $V(G)$. A graph which is not connected is said to be \emph{disconnected}. A \textit{connected component} of a graph $G$ is a maximal connected subgraph of $G$. For a graph $G$, we denote by $G-v$ the graph obtained from $G$ in which the vertex $v$ and all incident edges have been removed, and we denote by $G\setminus e$ the graph obtained from $G$ by removing the edge $e$. A \emph{cut vertex} of a connected graph $G$ is a vertex $v\in V(G)$ such that $G-v$ is disconnected. A \emph{bridge} of a connected graph $G$ is an edge $e\in E(G)$ for which $G \setminus e$ is disconnected.

A \textit{tree} is a connected graph that has no cycles. A \textit{forest} is a graph whose connected components are trees. A \textit{spanning tree} (resp. a \textit{spanning forest}) of a graph $G$ is a subgraph that is a tree (resp. a forest) and includes all of the vertices of $G$. A \textit{$k$-tree spanning forest} of $G$ is a spanning forest that consists of $k$ trees.

We also introduce some useful matrix and vector notation which will be used throughout. Let $G$ be a graph, and $H$ be a subgraph of $G$. If the vertices in $V(G)$ are labelled $v_1, \ldots, v_n$, we define $\mathbf{d}_G$ to be the column vector whose $i^\text{th}$ component is $\mathrm{deg}_G(v_i)$ for $1\leq i\leq n$. We use $\widehat{\mathbf{d}}_{H}$ to denote the vector $\widehat{\mathbf{d}}_{H}=(\hat{d}_i)_{v_i\in V(G)}$ where $\hat{d}_i=\mathrm{deg}_H(v_i)$ if $v_i\in V(H)$, and $\hat{d}_i=0$ if $v_i\in V(G)\backslash V(H)$. We denote by $\mathbf{1}_k$ the all-ones vector of length $k$, by $\mathbf{0}_k$ the all-zeros vector of length $k$. The subscript $k$ is omitted if the size is clear from the context. We denote by $\be_i$ the column vector whose $i^\text{th}$ entry is $1$ and zeros elsewhere. The size of $\be_i$ will be clear from the context.

While this work concerns itself mostly with combinatorial expressions for Markov chain quantities that can be derived in the case that our Markov chain is a random walk on a graph, we introduce the general concepts briefly here, as they are important in order to develop some intuition around the interpretation and meaning of Kemeny's constant. For further background, the interested reader is referred to \cite{KemenySnell}. 

A finite, discrete-time, time-homogeneous Markov chain is a stochastic process with a state space $\{s_1, \ldots, s_n\}$. At any given time the chain occupies a state $s_i$, and in discrete time-steps transitions to another state $s_j$ according to some prescribed transition probability $p_{ij}$. The probability transition matrix $P = [p_{ij}]$ determines the evolution of the system, in that the $(i,j)$ entry of $P^k$ is the probability of being in $s_j$ after $k$ time-steps, given that the chain starts in $s_i$. Under certain conditions on the Markov chain, the matrix $P^k$ converges to $\mathbf{1}\bm{\pi}^T$; that is, the probability distribution after $k$ steps converges to a stationary distribution vector $\bm{\pi}^T = \begin{bmatrix} \pi_1 & \pi_2 & \ldots & \pi_n \end{bmatrix}$, such that $\sum_i\pi_i = 1$. Thus $\pi_i$ can be interpreted as the long-term probability that the system occupies the state $s_i$. Note that $\bm{\pi}$ can be computed as a left Perron eigenvector for $P$, corresponding to the Perron value $1$. For short-term behaviour of the system, one can consider the \emph{mean first passage times}; for any pair of states $s_i, s_j$, the \emph{mean first passage time from $s_i$ to $s_j$}, denoted $m_{ij}$, is the expected number of time-steps elapsed until the chain reaches $s_j$ for the first time, given that it starts in $s_i$. These quantities can also be calculated via the transition matrix $P$, using the following formula:
\[m_{i,j} = \left\{\begin{array}{cc} \be_i^T (I-P_{(j)})^{-1}\mathbf{1}, & \mbox{if $i<j$};\\ \be_{i-1}^T(I-P_{(j)})^{-1}\mathbf{1}, & \mbox{if $i > j$},\end{array}\right.\]
where $P_{(j)}$ denotes the $j^{th}$ principal submatrix of $P$.

For a random walk on a graph $G$, the states of the Markov chain are the vertices of $G$, labelled $v_1, v_2, \ldots, v_n$. At any given time, the `random walker' occupies one of the vertices, and in each subsequent time-step, chooses a neighbour of the current vertex uniformly-at-random and moves there. Thus the transition probability $p_{ij}$ is given by 
\[p_{ij} = \left\{\begin{array}{cc} \frac{1}{\deg_G(v_i)}, & \mbox{if } v_i  \sim v_j;\\ 0, & \mbox{otherwise.}\end{array}\right.\]
Note that the stationary distribution vector for the simple random walk on a graph $G$ has $\pi_i = \frac{\deg_G(v_i)}{2|E(G)|}$, or $\bm{\pi} = \frac{1}{2m_G}\mathbf{d}_G$; that is, the long-term probability that a random walker finds themselves on $v_i$ is proportional to the vertex degree. 

Kemeny's constant is defined for an irreducible Markov chain by fixing an index $i$, and computing $\sum_{\substack{j=1\\j\neq i}}^n \pi_jm_{ij}$. As written, this can be interpreted as a weighted average of the mean first passage times from a fixed starting state. This is the astonishing quantity which is found to be independent of $i$; furthermore, since $\sum_i\pi_i = 1$, one can rewrite this as 
\[\kappa(P) = \sum_{i=1}^n \sum_{\substack{j=1\\j\neq i}}^n \pi_i m_{ij}\pi_j,\]
which provides the interpretation as the expected length of a random trip between states in the Markov chain, where both the starting and ending states are chosen at random, with respect to the stationary probability distribution.

Our final Markov chain parameter we define in the general setting is the \emph{accessibility index} of a state $s_j$ in an irreducible Markov chain. This is defined (see \cite{kirkland2016random}) as 
\[\alpha_j= \sum_{\substack{i=1\\i \neq j}}^n \pi_i m_{i,j}.\]
This appears at first glance to be very similar to Kemeny's constant; however, it is  a weighted average of the mean first passage times into a fixed destination state $s_j$, rather than from a fixed starting state. This quantity is not independent of the choice of index $j$, but rather defines a quantity to measure how easily a state $s_j$ is accessed from anywhere in the Markov chain (hence the name). Note that in the case of random walks on graphs, there is a related parameter known as the \emph{random walk centrality} of a vertex (see \cite{noh2004random}), which can be expressed as $\frac{1}{\alpha_j}$ (see \cite{kirkland2016random}). We note that $\sum_j \pi_j\alpha_j = \kappa(P)$.

For our work, we use a combinatorial expression for Kemeny's constant for a random walk on a connected and undirected graph which can be found in \cite{KirklandZeng}. We note that in \cite{KirklandZeng} (and also in \cite{PitmanTang}), there is a combinatorial expression given for $\kappa(P)$ for a general Markov chain using the all-minors matrix tree theorem; we refer the interested reader to these for extensions to weighted graphs or directed graphs, or arbitrary Markov chains with some interesting combinatorial structure in the transition matrix. 

In order to emphasize that we are dealing with random walks on connected and undirected graphs, given a connected graph $G$, we use $\kappa(G)$ to denote Kemeny's constant for the transition matrix of the random walk on $G$. We denote by $\tau_{G}$ the number of spanning trees of $G$, and by $\mathcal{F}_{G}(i;j)$ the set of $2$-tree spanning forests of $G$ such that one of the two trees contains a vertex $i$ of $G$, and the other has a vertex $j$ of $G$. This is occasionally referred to as a $2$-tree spanning forest \emph{separating} the vertices $i$ and $j$. Let  $f_{i,j}^{G}=|\mathcal{F}_{G}(i;j)|$, and define $F_{G}$ to be the matrix given by $F_{G}=[f_{i,j}^{G}]$. Then, $\kappa(G)$ is given by
\begin{equation}\label{formula:KemenyF}
\kappa(G)=\frac{\mathbf{d}_G^T F_G\mathbf{d}_G}{4m_G\tau_G}.
\end{equation}
Moreover, $\kappa(G)$ can be also expressed in terms of the \textit{effective resistance matrix} $R_G=[r_{i,j}^G]$ whose $(i,j)$-entry is defined as $r_{i,j}^G=(\be_i^T-\be_j^T)L^\dagger(\be_i-\be_j)$, where $L^\dagger$ is the Moore--Penrose inverse of the Laplacian matrix of $G$ (see \cite{bapat2010graphs}). The quantity $r_{i,j}^G$ is referred to as the \textit{effective resistance} or \emph{resistance distance} between vertices $i$ and $j$. It appears in \cite{chebotarev2020hitting} that $R_G=\frac{1}{\tau_G}F_G$. Hence, we also have
\begin{equation}\label{formula:KemenyR}
\kappa(G)=\frac{\mathbf{d}_G^T R_G\mathbf{d}_G}{4m_G}.
\end{equation}

This paper is organized as follows. In \Cref{Sec2}, we provide a formula of Kemeny's constant for graphs with a single bridge. With this formula, we examine how deleting the bridge and adding a bridge in a different place between the resulting components affects Kemeny's constant. In \Cref{Sec3}, we prove the main result of this article, deriving a formula for Kemeny's constant for graphs with multiple bridges. As done in \Cref{Sec2}, we investigate how to optimize Kemeny's constant by the configuration of the bridges. Finally, we close this paper by outlining several potential applications in \Cref{Sec4}.

\section{Kemeny's constant for graphs with a cut vertex or a bridge}\label{Sec2}

The combinatorial building-blocks for the formula for $\kappa(G)$ in \eqref{formula:KemenyF} are the degree vector $\mathbf{d}_G$, the number of edges $m_G$, the number of spanning trees $\tau_G$, and the matrix $F_G$ of numbers of 2-tree spanning forests of $G$ separating $i$ and $j$. The following proposition gives expressions for these building-blocks in the case that $G$ is constructed by connecting two graphs with a bridge. The proof for this proposition is an adaptation of the proof of the result \cite[Prop 3.1]{kim2022families} describing similar quantities for a graph with a cut vertex. It is included here for completeness.

For a graph $G$, we use $\mathbf{f}_G^j$ to denote the $j^{\text{th}}$ column of $F_{G}$, or the column of $F_G$ corresponding to vertex $j$.

\begin{proposition}\label{Proposition:dRd bridge}
	Let $G$ be a connected graph on $n$ vertices, and let $e$ be a bridge in $G$. Suppose that $G_1$ and $G_2$ are the components of $G\backslash e$. Let $e=v_1\sim v_2$ where $v_1\in V(G_1)$ and $v_2\in V(G_2)$. Then, labelling the vertices of $G$ in order of $V(G_1)$ and $V(G_2)$, we have:
	\begin{align*}
	\mathbf{d}_{G}&=\widehat{\mathbf{d}}_{G_1}+\widehat{\mathbf{d}}_{G_2}+\be_{v_1}+\be_{v_2},\\m_{G}&=m_{G_1}+m_{G_2}+1,\\\tau_{G}&=\tau_{G_1}\tau_{G_2},\\
	F_{G}&=\left[\begin{array}{c|c}
	\tau_{G_2}F_{G_1} & \tau_{G_2}\mathbf{f}_{G_1}^{v_1}\mathbf{1}^T\hspace*{-0.1cm}+\tau_{G_1}\mathbf{1}(\mathbf{f}_{G_{2}}^{v_2})_{\phantom{1_{1_1}}}^T\hspace*{-0.3cm}+\tau_{G_1}\tau_{G_{2}}J \\\hline
	\tau_{G_1}\mathbf{f}_{G_2}^{v_2}\mathbf{1}^T\hspace*{-0.1cm}+\tau_{G_2}\mathbf{1}(\mathbf{f}_{G_1}^{v_1})^{T^{\phantom{T^a}}}\hspace*{-0.4cm}+\tau_{G_1}\tau_{G_2}J & \tau_{G_1}F_{G_2}
	\end{array}\right].
	\end{align*}
	Moreover, this implies that
	\begin{align*}
	\mathbf{d}_G^TF_G\mathbf{d}_G=&\tau_{G_2}\mathbf{d}_{G_1}^TF_{G_1}\mathbf{d}_{G_1}+\tau_{G_1}\mathbf{d}_{G_2}^TF_{G_2}\mathbf{d}_{G_2}+4\tau_{G_2}(m_{G_2}+1)\mathbf{d}_{G_1}^T\mathbf{f}_{G_1}^{v_1}\\
	&+4\tau_{G_1}(m_{G_1}+1)\mathbf{d}_{G_2}^T\mathbf{f}_{G_2}^{v_2}+2\tau_{G_1}\tau_{G_2}(2m_{G_1}+1)(2m_{G_2}+1);
	\end{align*}
	and equivalently,
	\begin{align}\label{1;ExpressionR}
	\mathbf{d}_G^TR_{G}\mathbf{d}_G=&\mathbf{d}_{G_1}^TR_{G_1}\mathbf{d}_{G_1}+\mathbf{d}_{G_2}^TR_{G_2}\mathbf{d}_{G_2}+4(m_{G_2}+1)\mathbf{d}_{G_1}^TR_{G_1}\be_{v_1}\\\nonumber
	&+4(m_{G_1}+1)\mathbf{d}_{G_2}^TR_{G_2}\be_{v_2}+2(2m_{G_1}+1)(2m_{G_2}+1).
	\end{align}
\end{proposition}
\begin{proof}
	The results for $\mathbf{d}_{G}$ and $m_G$ are straightforward. Since any spanning tree of $G$ contains the bridge $e$, we have $\tau_G=\tau_{G_1}\tau_{G_2}$. Note that $F_G$ is symmetric. We consider two cases for the structure of $F_G$: \begin{enumerate*}[label=(\roman*)]
		\item\label{1;cond1} vertices $i$ and $j$ both are in either $V(G_1)$ or $V(G_2)$;
		\item\label{1;cond2} one of $i$ and $j$ is in $V(G_1)$, and the other is in $V(G_2)$.
	\end{enumerate*}
	Consider \ref{1;cond1}. Let $i,j\in V(G_1)$. Then, any $2$-tree spanning forest in $\mathcal{F}_{G}(i;j)$ can be constructed from a $2$-tree spanning forest in $\mathcal{F}_{G_1}(i;j)$ and a spanning tree of $G_2$ by joining $v_1$ and $v_2$, and vice versa. Thus, $|\mathcal{F}_{G}(i;j)|=\tau_{G_2}|\mathcal{F}_{G_1}(i;j)|$. In a similar way, we can find $|\mathcal{F}_{G}(i;j)|=\tau_{G_1}|\mathcal{F}_{G_2}(i;j)|$ for $i,j\in V(G_2)$. 
	
	For the case \ref{1;cond2}, we assume without loss of generality that $i\in V(G_1)$ and $j\in V(G_2)$. Then, $\mathcal{F}_{G}(i;j)$ is the union of three disjoint subsets $X_1$, $X_2$, and $X_3$, where $X_1$ (resp. $X_2$) is the set of $2$-tree spanning forests such that the subtree having the vertex $i$ (resp. $j$) contains $v_2$ (resp. $v_1$); and $X_3$ is the set of $2$-tree spanning forests such that the subtrees having $i$ and $j$ contain $v_1$ and $v_2$, respectively. Pick any forest in $X_1$. Then, it can be constructed from a forest in $\mathcal{F}_{G_2}(v_2;j)$ and a spanning tree of $G_1$, by joining $v_1$ and $v_2$, and vice versa. So, we have $|X_1|=\tau_{G_1}|\mathcal{F}_{G_2}(v_2;j)|$. Similarly, $|X_2|=\tau_{G_2}|\mathcal{F}_{G_1}(i;v_1)|$. Consider $X_3$. Since any $2$-tree spanning forest in $X_3$ does not contain $e$, we have $|X_3|=\tau_{G_1}\tau_{G_2}$. Therefore, our desired structure of $F_G$ is established.
	
	For computation of $\mathbf{d}_G^TF_G\mathbf{d}_G$, one can check that
	\begin{align*}
	(\be_{v_1}\hspace*{-0.12cm}+\be_{v_2})^T F_G (\be_{v_1}\hspace*{-0.12cm}+\be_{v_2})&=2\be_{v_1}^TF_G\be_{v_2}=2\tau_{G_1}\tau_{G_2}\\
	[\mathbf{d}_{G_1}^T\;\mathbf{0}^T]F_G(\be_{v_1}\hspace*{-0.12cm}+\be_{v_2})&=2\tau_{G_2}\mathbf{d}_{G_1}^T\mathbf{f}_{G_1}^{v_1}\hspace*{-0.1cm}+2\tau_{G_1}\tau_{G_2}m_{G_1}\\
	[\mathbf{0}^T\;\mathbf{d}_{G_2}^T]F_G(\be_{v_1}\hspace*{-0.12cm}+\be_{v_2})&=2\tau_{G_1}\mathbf{d}_{G_2}^T\mathbf{f}_{G_2}^{v_2}\hspace*{-0.1cm}+2\tau_{G_1}\tau_{G_2}m_{G_2}\\
	[\mathbf{d}_{G_1}^T\;\mathbf{0}^T]F_G\begin{bmatrix}
	\mathbf{0}\\
	\mathbf{d}_{G_2}
	\end{bmatrix}&=2m_{G_2}\tau_{G_2}\mathbf{d}_{G_1}^T\mathbf{f}_{G_1}^{v_1}\hspace*{-0.1cm}+2m_{G_1}\tau_{G_1}\mathbf{d}_{G_2}^T\mathbf{f}_{G_2}^{v_2}\hspace*{-0.1cm}+4\tau_{G_1}\tau_{G_2}m_{G_1}m_{G_2}.
	\end{align*}
	Then, the remaining conclusion follows. 
\end{proof}

Using these building blocks, we can express Kemeny's constant for a graph with a bridge.

\begin{proposition}\label{Prop1:bridge between 2}
	Let $G$ be a connected graph, and let $e$ be a bridge in $G$. Suppose that $G_1$ and $G_2$ are the components of $G\backslash e$. Let $e=v_1\sim v_2$, $v_1\in V(G_1)$, and $v_2\in V(G_2)$. Then,
	\begin{align*}
	\kappa(G)=&\frac{m_{G_1}}{m_{G}}\kappa(G_1)+\frac{m_{G_2}}{m_G}\kappa(G_2)+\frac{m_{G_2}+1}{m_G}\mathbf{d}_{G_1}^TR_{G_1}\be_{v_1}\\
	&+\frac{m_{G_1}+1}{m_G}\mathbf{d}_{G_2}^TR_{G_2}\be_{v_2}+\frac{(2m_{G_1}+1)(2m_{G_2}+1)}{2m_G}.
	\end{align*}
\end{proposition}
\begin{proof}
	Recall that $\kappa(G)=\frac{\mathbf{d}_G^T R_G\mathbf{d}_G}{4m_G}$. Dividing both sides of \eqref{1;ExpressionR} by $4m_G$, one can establish the desired result.
\end{proof}

Given two graphs $G_1$ and $G_2$, it is interesting to consider the value of Kemeny's constant of the resultant graph $G$ when a bridge $v_1\sim v_2$ has been added for some $v_1 \in V(G_1)$ and some $v_2 \in V(G_2)$, and to explore how the choice of $v_1$ and $v_2$ can affect the value of $\kappa(G)$. In particular, the range of possible values of $\kappa(G)$ hinges on the possible values of the quantities $\mathbf{d}_{G_1}^TR_{G_1}\be_{v_1}$ and $\mathbf{d}_{G_2}^TR_{G_2}\be_{v_2}$. This quantity is more formally defined in \cite{ciardo2020kemeny,faught20221} as follows:

\begin{definition}
	Let $G$ be a connected graph, and $v$ be a vertex of $G$. The \textit{moment} $\mu_G(v)$ of $v$ in $G$ is defined as
	\begin{align*}
	\mu_G(v)=\mathbf{d}_{G}^TR_{G}\be_{v}.
	\end{align*}
\end{definition}

The following lemma indicates how the moment $\mu_G(v)$ can be expressed in terms of $\kappa(G)$ and the accessibility index $\alpha_G(v)$ in $G$, which is defined in \cite{kirkland2016random} as the $v^\text{th}$ entry of $w^TM-\mathbf{1}^T$ where $w$ and $M$ are the stationary vector and mean first passage time matrix, respectively, for the random walk on $G$.

\begin{lemma}\label{lemma:acce moment kemeny}
	Let $G$ be a connected graph. Then 
	\[\mu_G(v) = \alpha_G(v) + \kappa(G).\]
\end{lemma}
\begin{proof}
	Note that 
	\[\mathbf{d}_G^T R_G e_v = \sum_{i=1}^n \deg_G(i)r_{i,v}^G.\]
	It is shown in \cite{chandra1996electrical} that the effective resistance between two vertices $i$ and $j$ satisfies
	\[2m_G\cdot r_{i,j}^G = m_{i,j} + m_{j,i},\]
	where $m_{i,j}$ is the mean first passage time from $i$ to $j$ (note that $m_{i,j}+m_{j,i}$ is sometimes referred to as the \emph{commute time} of vertices $i$ and $j$). Recall that $\bm{\pi}=\frac{1}{2m}\mathbf{d}_G$ is the stationary vector of a random walk on $G$. Thus
	\begin{eqnarray*}
		\sum_{i=1}^n \deg_G(i)r_{i,j}^G & = & \sum_{\substack{i=1\\i\neq j}}^n \frac{\deg_G(i)}{2m_G} m_{i,j} + \sum_{\substack{i=1\\i\neq j}}^n \frac{\deg_G(i)}{2m_G} m_{ji} \\
		& = & \sum_{\substack{i=1\\i\neq j}}^n \pi_i m_{i,j} + \sum_{\substack{i=1\\i\neq j}}^n \pi_i m_{ji} \\
		& = & \alpha_{G}(j) + \kappa(G).
	\end{eqnarray*}
\end{proof}

Combining this observation with Proposition \ref{Prop1:bridge between 2}, we achieve the following alternative expression for Kemeny's constant of a graph with a bridge, expressed in terms of the values for Kemeny's constant of the connected components $G_1$ and $G_2$ upon removal of the bridge, of the number of edges in each component $m_{G_1}$ and $m_{G_2}$, and of the accessibility index of each vertex incident with the bridge in each component.

\begin{theorem}\label{kemeny;two components}
	Let $G$ be a connected graph with a bridge $e = \{v_1, v_2\}$, and let $G_1$ and $G_2$ be the connected components of $G\setminus e$, with $v_1 \in V(G_1)$ and $v_2\in V(G_2)$. Then 
	\begin{align*}
	\kappa(G) = & \; \kappa(G_1) + \kappa(G_2) + \frac{m_{G_2}+1}{m_G}\alpha_{G_1}(v_1) + \frac{m_{G_1}+1}{m_G} \alpha_{G_2}(v_2)\\ & \quad +  \frac{(2m_{G_1}+1)(2m_{G_2}+1)}{2m_G}.
	\end{align*}
\end{theorem}

\subsection{Optimization of Kemeny's constant for graphs with a bridge}\label{Subsec:optimization1}

In this section, we discuss how \Cref{kemeny;two components} provides insight on the range of possible values of Kemeny's constant for a graph $G$ created by joining two connected graphs $G_1$ and $G_2$ with a bridge $v_1\sim v_2$, where $v_1 \in V(G_1)$ and $v_2 \in V(G_2)$. In particular, in this subsection we consider how to maximize and minimize Kemeny's constant for such a graph; from \Cref{kemeny;two components}, it is equivalent to maximizing/minimizing $\alpha_{G_1}(v_1)$ and $\alpha_{G_2}(v_2)$. 

\begin{example}\label{Example_Ema}
	Let $G_1$ and $G_2$ be the following graphs whose vertices are labelled by their accessibility indices:
	
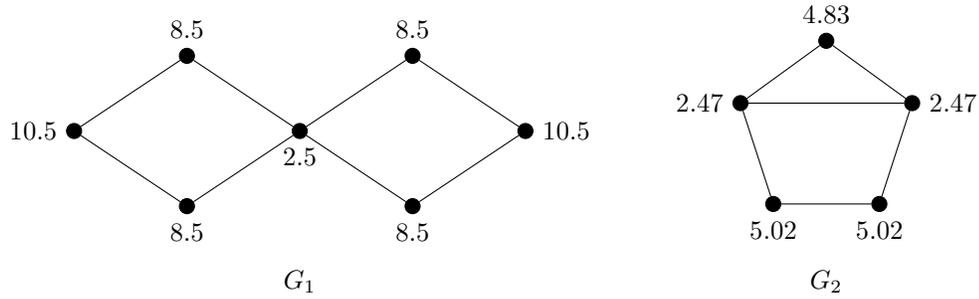
\begin{center}
		\begin{tikzpicture}
		\tikzset{enclosed/.style={draw, circle, inner sep=0pt, minimum size=.20cm, fill=black}}
		
		\node[enclosed, label={below, yshift=0cm: $2.5$}] (v_1) at (0,0) {};
		\node[enclosed, label={below, yshift=0cm: $8.5$}] (v_2) at (-1.5,-1) {};
		\node[enclosed, label={above, yshift=0cm: $8.5$}] (v_3) at (-1.5,1) {};
		\node[enclosed, label={left, xshift=0cm: $10.5$}] (v_4) at (-3,0) {};
		\node[enclosed, label={below, yshift=0cm: $8.5$}] (v_5) at (1.5,-1) {};
		\node[enclosed, label={above, yshift=0cm: $8.5$}] (v_6) at (1.5,1) {};
		\node[enclosed, label={right, xshift=0cm: $10.5$}] (v_7) at (3,0) {};
		
		\draw (v_1) -- (v_2);
		\draw (v_1) -- (v_3);
		\draw (v_2) -- (v_4);
		\draw (v_3) -- (v_4);
		
		\draw (v_1) -- (v_5);
		\draw (v_1) -- (v_6);
		\draw (v_5) -- (v_7);
		\draw (v_6) -- (v_7);
		
		\node[] at (0,-2) {$G_1$};
		
		\node[enclosed, label={right,xshift=0cm: $2.47$}] (1) at ({7+1.2*cos(18+72*(0))}, {1.2*sin(18+72*(0))}) {};
		\node[enclosed, label={above,yshift=0cm: $4.83$}] (2) at ({7+1.2*cos(18+72*(1))}, {1.2*sin(18+72*(1))}) {};
		\node[enclosed, label={left,xshift=0cm: $2.47$}] (3) at ({7+1.2*cos(18+72*(2))}, {1.2*sin(18+72*(2))}) {};
		\node[enclosed, label={below,yshift=0cm: $5.02$}] (4) at ({7+1.2*cos(18+72*(3))}, {1.2*sin(18+72*(3))}) {};
		\node[enclosed, label={below,yshift=0cm: $5.02$}] (5) at ({7+1.2*cos(18+72*(4))}, {1.2*sin(18+72*(4))}) {};
		\draw (1)--(2)--(3)--(4)--(5)--(1);
		\draw (1)--(3);
		\node[] at (7, -2) {$G_2$};
		\end{tikzpicture}
	\end{center}
	Suppose that $G$ is formed from $G_1$ and $G_2$ by connecting a vertex of $G_1$ and a vertex from $G_2$ with an edge. Note that $\kappa(G_1) = 7.5$, $\kappa(G_2)=3.71$, $m_{G_1} = 8$, $m_{G_2}=6$, and $m_G$ will be equal to 15. These quantities impose restrictions on the range of values of $\kappa(G)$, and thus the minimum will occur when the vertices of minimum accessibility index in each graph are joined by a bridge, and the maximum will occur when the vertices of maximum accessibility index in each graph are joined. Thus
	\[20.687 = \kappa(G_{min}) \leq \kappa(G) \leq \kappa(G_{max})= 25.947,\]
	where $G_{min}$ and $G_{max}$ are as shown.
	
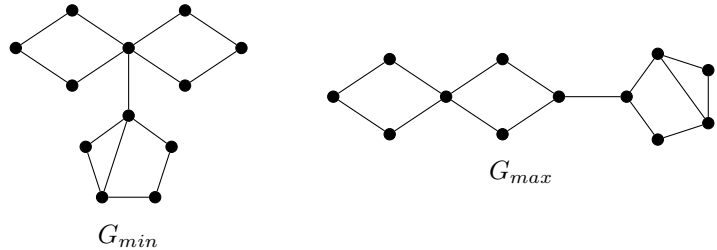
\begin{center}
		\begin{tikzpicture}[scale=0.5]
		\tikzset{enclosed/.style={draw, circle, inner sep=0pt, minimum size=.150cm, fill=black}}
		
		\node[enclosed] (v_1) at (0,0) {};
		\node[enclosed] (v_2) at (-1.5,-1) {};
		\node[enclosed] (v_3) at (-1.5,1) {};
		\node[enclosed] (v_4) at (-3,0) {};
		\node[enclosed] (v_5) at (1.5,-1) {};
		\node[enclosed] (v_6) at (1.5,1) {};
		\node[enclosed] (v_7) at (3,0) {};
		
		\draw (v_1) -- (v_2);
		\draw (v_1) -- (v_3);
		\draw (v_2) -- (v_4);
		\draw (v_3) -- (v_4);
		
		\draw (v_1) -- (v_5);
		\draw (v_1) -- (v_6);
		\draw (v_5) -- (v_7);
		\draw (v_6) -- (v_7);
		
		\foreach \x in {1,...,5} \node[enclosed] (\x) at ({1.2*cos(18+72*(\x-1))}, {-3+1.2*sin(18+72*(\x-1))}) {};
		\draw (1)--(2)--(3)--(4)--(5)--(1);
		\draw (2)--(4);
		\draw (2)--(v_1);
		
		\node[] at (0,-5) {$G_{min}$};
		\end{tikzpicture}\hspace{30pt}\begin{tikzpicture}[scale=0.5]
		\tikzset{enclosed/.style={draw, circle, inner sep=0pt, minimum size=.150cm, fill=black}}
		
		\node[enclosed] (v_1) at (0,0) {};
		\node[enclosed] (v_2) at (-1.5,-1) {};
		\node[enclosed] (v_3) at (-1.5,1) {};
		\node[enclosed] (v_4) at (-3,0) {};
		\node[enclosed] (v_5) at (1.5,-1) {};
		\node[enclosed] (v_6) at (1.5,1) {};
		\node[enclosed] (v_7) at (3,0) {};
		
		\draw (v_1) -- (v_2);
		\draw (v_1) -- (v_3);
		\draw (v_2) -- (v_4);
		\draw (v_3) -- (v_4);
		
		\draw (v_1) -- (v_5);
		\draw (v_1) -- (v_6);
		\draw (v_5) -- (v_7);
		\draw (v_6) -- (v_7);
		
		\foreach \x in {1,...,5} \node[enclosed] (\x) at ({6+1.2*cos(36+72*(\x-1))}, {1.2*sin(36+72*(\x-1))}) {};
		\draw (1)--(2)--(3)--(4)--(5)--(1);
		\draw (2)--(5);
		\draw (3)--(v_7);
		\node[] at (2,-2) {$G_{max}$};
		\node[] at (2, -4) {\,};
		\end{tikzpicture}
	\end{center}
\end{example}

Now we shall show an interesting result that if $G_1$ and $G_2$ are chosen as trees in \Cref{kemeny;two components}, then the minimum of $\kappa(G)$ is attained when $v_1$ and $v_2$ are \textit{centroids} of $G_1$ and $G_2$, respectively. In \cite{kirkland2016random}, the author examines characterizations of the minimum and maximum of accessibility index in a tree; that is, the minimum is attained at either a unique vertex or two adjacent vertices, while the maximum only at a pendent vertex. Here we show at which vertex we can attain the minimum value of the accessibility indices. A vertex $v$ of a tree $\cT$ on $n$ vertices is called a \textit{centroid} if each subtree in $\cT-v$ contains at most $\lfloor\frac{n}{2}\rfloor$ vertices. It is well-known that either $\cT$ contains a unique centroid $v$ so that each subtree in $\cT-v$ contains less than $\frac{n}{2}$, or $\cT$ has exactly two adjacent centroids $v_1$ and $v_2$, in which case $n$ is even, and two components in $\cT\backslash v_1\sim v_2$ contain $\frac{n}{2}$ vertices, respectively.

\begin{proposition}\label{centroid minimum}
	Let $\cT$ be a tree. Then, $\min_{u \in V(\cT)} \alpha_\cT(u)$ is attained at $v$ if and only if $v$ is a centroid.
\end{proposition}
\begin{proof}
	It is found in \cite[Lemma 9.7]{bapat2010graphs} that $\mathbf{d}_\cT^T F_\cT=\mathbf{1}^T(2F_\cT-(n-1)I)$. Thus, the minimum value of the accessibility index is attained at a vertex corresponding to the minimum entry in $\mathbf{1}^T F_\cT$.
	
	For any two adjacent vertices $v_1,v_2\in V(\cT)$, let $\cT_{v_1}$ and $\cT_{v_2}$ be the components in $\cT\backslash v_1\sim v_2$ that contain $v_1$ and $v_2$, respectively. We note that the $(i,j)$-entry in $F_\cT$ is the distance between vertex $i$ and $j$ in $\cT$. Then,
	\begin{align*}
	\mathbf{1}^T F_{\cT}\be_{v_2}&= \sum_{x\in V(\cT)} \mathrm{dist}_{\cT}(v_2, x) \\
	&=\sum_{x\in V(\cT_{v_1})}\left(\mathrm{dist}_\cT(v_1,x)+1\right)+\sum_{x\in V(\cT_{v_2})}(\mathrm{dist}_\cT(v_1,x)-1)\\
	&=\sum_{x\in V(\cT)}\mathrm{dist}_\cT(v_1,x) + |V(\cT_{v_1})|-|V(\cT_{v_2})|.
	\end{align*}
	Hence,
	\begin{align}\label{temp:identity}
	\mathbf{1}^T F_{\cT}\be_{v_1}-\mathbf{1}^T F_{\cT}\be_{v_2}=|V(\cT_{v_2})|-|V(\cT_{v_1})|.
	\end{align}
	
	Let $v$ be a centroid of $G$. We consider a vertex $w$ of $\cT$ adjacent to $v$. Then, $|V(\cT_{w})|\leq \frac{n}{2}$. If $|V(\cT_{w})|=\frac{n}{2}$, then $w$ is the other centroid and using \eqref{temp:identity} we have $\mathbf{1}^T F_{\cT}\be_{v}=\mathbf{1}^T F_{\cT}\be_{w}$. If $|V(\cT_{w})|<\frac{n}{2}$, then it follows from \eqref{temp:identity} that $\mathbf{1}^T F_{\cT}\be_{v}<\mathbf{1}^T F_{\cT}\be_{w}$. Now we let $v_1, v_2$ be any two adjacent vertices distinct from $v$. Suppose without loss of generality that $1\leq \mathrm{dist}_\cT(v,v_1)<\mathrm{dist}_\cT(v,v_2)$. We claim that $\mathbf{1}^T F_{\cT}\be_{v_1}<\mathbf{1}^T F_{\cT}\be_{v_2}$. Assume to the contrary $|V(\cT_{v_2})|-|V(\cT_{v_1})|\geq 0$. Since $|V(\cT_{v_2})|+|V(\cT_{v_1})|=n$, we have $|V(\cT_{v_2})|\geq\frac{n}{2}$, and so the subtree in $\cT-v$ containing $v_2$ has at least $\frac{n}{2}+1$ vertices, which contradicts that $v$ is a centroid. Hence, $|V(\cT_{v_2})|-|V(\cT_{v_1})|<0$, and $\mathbf{1}^T F_{\cT}\be_{v_1}<\mathbf{1}^T F_{\cT}\be_{v_2}$. Therefore, we obtain our desired result.
\end{proof}

\begin{corollary}\label{minimum of two components}
	Let $\cT_1$ and $\cT_2$ be trees, and let $G$ be a graph formed by adding an edge $v_1\sim v_2$, where $v_1 \in V(\cT_1)$ and $v_2\in V(\cT_2)$. Then, $\kappa(G)$ is minimized if and only if $v_1$ and $v_2$ are centroids of $\cT_1$ and $\cT_2$, respectively.
\end{corollary}
\begin{proof}
	The conclusion follows from \Cref{kemeny;two components} and \Cref{centroid minimum}.
\end{proof}

\section{Kemeny's constant for a chain of connected graphs with respect to a tree}\label{Sec3}


\begin{definition}\label{Def:a Chain}
	Let $\cT$ be a tree on $k$ vertices where $V(\cT)=\{1,\dots,k\}$. Let $G_1,\dots,G_k$ be connected graphs. Let $G$ be a graph constructed as follows: the vertices $1,\dots,k$ are replaced by the graphs $G_1,\dots,G_k$, respectively; and if $i\sim j$ is an edge of $\cT$, then some vertex $v_i\in V(G_i)$ is chosen, and some vertex $v_j\in V(G_j)$ is chosen, and the two vertices are joined with an edge so that $v_i\sim v_j$ is a bridge in $G$. Then $G$ is said to be a \textit{chain of $G_1,\dots,G_k$ with respect to $\cT$}. We denote by $\cB_G$ the set of the $(k-1)$ bridges, used in the construction of $G$, that correspond to the edges of $\cT$.
\end{definition}

\begin{figure}[h!]
	\centering
	\begin{tikzpicture}[scale=.6]
	\newcommand\Square{+(-1,-1) rectangle +(1,1)}
	
	\tikzset{enclosed/.style={draw, circle, inner sep=0pt, minimum size=.20cm, fill=black}}
	\def \radius {3cm}
	\node[enclosed, label=below : $1$] (ustar) at (360:0mm) {};
	\foreach \i [count=\ni from 0] in {4,3,2,k,k-1}{
		\node[enclosed, label=above: $\i$] at ({210-\ni*60}:\radius) (u\ni) {};
		\draw (ustar)--(u\ni);
	}
	
	\draw[thick, loosely dotted] (240:\radius) arc[start angle=240, end angle=300, radius=\radius];
	\node[] at (0,-5) {$\mathcal{T}$};
	\node[label=$ $] (p) at (1.5,-5) {};
	
	\begin{scope}[xshift=10cm]
	
	\def \radius {2.5cm}
	
	\def \margin {8} 
	
	\node[draw, circle] at (360:0mm) (ustar) {$G_1$};
	\foreach \i [count=\ni from 0] in {4,3,2,k}{
		\node[draw, circle] at ({210-\ni*60}:\radius) (u\ni) {$G_{\i}$};
		\draw (ustar)--(u\ni);
	}
	\node[outer sep=-1.5] at (-30:\radius) (u5) {$G_{k-1}$};
	\path[draw=black] (u5) circle[radius=0.8];
	\draw (ustar)--(u5);
	
	\draw[thick, loosely dotted] (240:\radius) arc[start angle=240, end angle=300, radius=\radius];
	\node[] at (0,-5) {$G$};
	\node[label=$ $] (p') at (-1.5,-5) {};
	\end{scope}
	
	\draw[->] (p) edge[bend right] (p');
	\end{tikzpicture}
	\caption{An example for illustration of Definition \ref{Def:a Chain}, which is a chain $G$ of connected graphs $G_1,\dots,G_k$ with respect to a star $\cT$ on $k$ vertices.}\label{Figure:Chain of Graphs}
\end{figure}
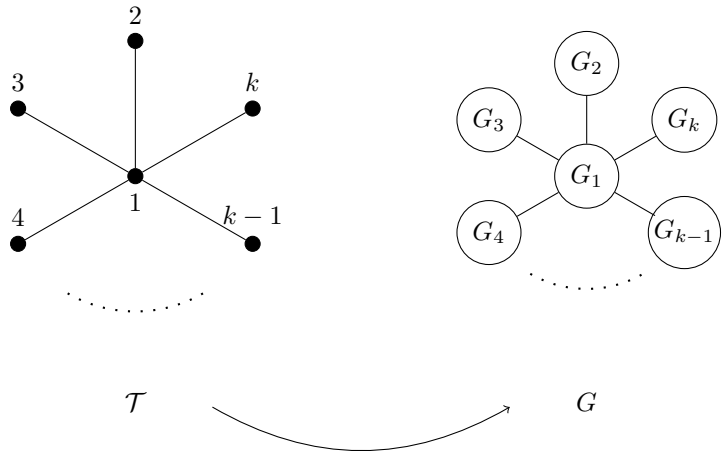

\begin{remark}\label{Remark:vertices for T and G}
	When we deal with a chain $G$ of graphs with respect to a tree $\cT$ on $k$ vertices, we write the vertex set of $\cT$ as $\{1,\dots,k\}$ and the vertices of $G$ are labelled by letters other than $1,\dots,k$. 
\end{remark}

In this section, we provide a formula for Kemeny's constant for a chain $G$ of connected graphs $G_1,\dots,G_k$ with respect to a tree. The formula in Theorem \ref{Theorem:MainResult} is proved by induction using Proposition \ref{Prop1:bridge between 2} and a lemma we prove in this section, Lemma \ref{Lemma2:dRe}. 
Since the formula for $\kappa(G)$ is based on resistance distances, which are in turn based on $2$-tree spanning forests, we will need to consider how such forests can be constructed from spanning trees and forests of $G_1, G_2, \ldots, G_k$, taking into account the structure of $\cT$. This is done in Proposition~\ref{Proposition:entries of f in chain}. For ease of exposition, we spend some time introducing and elaborating on the notation we use in this section, and list some as Observations for ease of reference later.

Let $\cT$ be a tree on $k$ vertices, and let $G$ be a chain of connected graphs $G_1,\dots,G_k$ with respect to $\cT$. For $v,w\in V(G)$, there exist $i, j\in\{1,\dots,k\}$ such that $v\in V(G_i)$ and $w\in V(G_j)$. Let $l=\mathrm{dist}_\cT(i,j)+1$, and suppose $l\geq 2$. Since $\cT$ is a tree, there is a unique path $(i_1,\dots,i_l)$ from $i$ to $j$ in $\cT$, that is, $i_m\sim i_{m+1}$ for $m=1,\dots,l-1$. (Here, $i_1=i$ and $i_l=j$.) Moreover, letting $v_1=v$ and $w_l=w$, there are exactly $l-1$ bridges in $\cB_G$ such that $w_{m}\sim v_{m+1}$ is the bridge between $G_{i_m}$ and $G_{i_{m+1}}$ for $m=1,\dots,l-1$. Hence, given $v,w\in V(G)$, we may define a set consisting of pairs of vertices,
\begin{align*}
P_{v,w}=\{(v_i,w_i)\,|\,i=1,\dots,l\}.
\end{align*}
We note that $v_1$ may be the same as $w_1$, so for this case, $(v_1,v_1)\in P_{v,w}$; similarly, $v_l$ may be the same as $w_l$. We further note that $P_{w_1,w}=\left(P_{v,w}\backslash\{(v,w_1)\}\right)\cup \{(w_1,w_1)\}$. For the case $l=1$, we have $i=j$, and $v, w$ are both in $V(G_i)$. In this case, we let $P_{v,w}=\{(v,w)\}$. The purpose of defining this set is for indexing certain summations in the results which follow; for example, in Proposition \ref{Proposition:entries of f in chain}, the number of $2$-tree spanning forest separating $v$ and $w$ can be calculated in terms of the number of $2$-tree spanning forests which separate pairs of vertices in $P_{v, w}$.

\begin{figure}[h!]
	\centering
	\begin{tikzpicture}
	\tikzset{enclosed/.style={draw, circle, inner sep=0pt, minimum size=.065cm, fill=black}}
	\node[enclosed, label={below : $v$}]  at (-4.8,0) {};
	\node[enclosed, xshift=0.6cm, label={below, xshift=.2cm : $w_1$}] (w1) at (-4.8,0) {};
	\node[enclosed, xshift=-0.6cm, label={below, xshift=-.2cm : $v_2$}] (v2) at (-2.4,0) {};
	\node[enclosed, xshift=0.6cm, label={below, xshift=.2cm : $w_2$}] (w2) at (-2.4,0) {};
	\node[] (v3) at (-.5,0) {}; 
	
	\node[] (wl-2) at (.5,0) {};
	\node[enclosed, xshift=-0.6cm, label={below, xshift=-.25cm : $v_{l-1}$}] (vl-1) at (2.4,0) {};
	\node[enclosed, xshift=0.6cm, label={below, xshift=.3cm : $w_{l-1}$}] (wl-1) at (2.4,0) {};
	\node[enclosed, xshift=-0.6cm, label={below, xshift=-.2cm : $v_l$}] (vl) at (4.8,0) {};
	\node[enclosed, label={below : $w$}]  at (4.8,0) {};

	\draw (-4.8,0) circle (0.6cm);
	\draw (-2.4,0) circle (0.6cm);
	\draw[thick, loosely dotted] (-.5,0)--(.5,0);
	\draw (2.4,0) circle (0.6cm);
	\draw (4.8,0) circle (0.6cm);
	
	\draw (w1) -- (v2); \draw (w2)--(v3);
	\draw (wl-2) -- (vl-1); \draw (wl-1)--(vl);
	
	\node[] at (-4.8,-1) {$G_{i_1}$};
	\node[] at (-2.4,-1) {$G_{i_2}$};
	\node[] at (2.4,-1) {$G_{i_{l-1}}$};
	\node[] at (4.8,-1) {$G_{i_l}$};
	\end{tikzpicture}
	\caption{Illustration of the definition of the set $P_{v,w}$. In this figure, we have $P_{v, w} = \{(v, w_1), (v_2, w_2), \ldots, (v_{l-1}, w_{l-1}), (v_l, w)\}$.}\label{Figure:P{v,w}}
\end{figure}
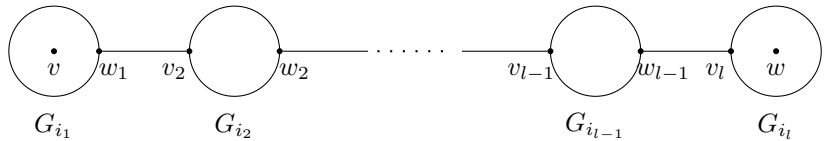

Given any vertex $z$ in $G$, there exists a unique $j\in\{1,\dots,k\}$ such that $z\in V(G_j)$; and so we shall use $G_z$ to denote the graph $G_j$ (as stated in Remark \ref{Remark:vertices for T and G}, $z$ is not a number between $1$ and $k$). We note that for any pair $(x,y)$ in $P_{v,w}$, we have $G_x=G_y$. 

\begin{observation}\label{Observation2}
	Let $\cT$ be a tree on $k$ vertices, and let $G$ be a chain of connected graphs $G_1,\dots,G_k$ with respect to $\cT$. Let $v,w\in V(G)$. Then, $v\in V(G_{k_1})$ and $w\in V(G_{k_2})$ for some $k_1,k_2\in\{1,\dots,k\}$. Consider the set $P_{v,w}$, and pick any pair $(x,y)\in P_{v,w}$. Then there exists a unique $\ell\in\{1,\dots,k\}$ such that $G_\ell=G_x$. Considering how $P_{v,w}$ is defined, $\ell$ is a vertex of $\cT$ that lies on the path from $i$ to $j$ in $\cT$. Therefore, elements in $P_{v,w}$ are in one-to-one correspondence with vertices on the path from $i$ to $j$ in $\cT$.
\end{observation}

The next proposition enumerates the number of $2$-tree spanning forests separating any two vertices $v$ and $w$ in $G$, and uses this expression to determine the resistance distances, and the quantity $\widehat{\mathbf{d}}_{G_i}^TR_{G}\be_{w}$. Note that the sizes of the vectors $\be_i$ are ambiguous but they can be determined from the context. Furthermore, the vector $\be_w$ denotes the standard basis vector with a 1 in the position corresponding to the vertex $w$.

\begin{proposition}\label{Proposition:entries of f in chain}
	Given a tree $\cT$ on $k$ vertices and connected graphs $G_1,\dots,G_k$, let $G$ be a chain of $G_1,\dots,G_k$ with respect to $\cT$. Consider $v,w\in V(G)$, and that $v\in V(G_i)$ and $w\in V(G_j)$ for some $i,j\in\{1,\dots,k\}$. Let $w_1$ be the vertex of $G_i$ such that $(v,w_1)\in P_{v,w}$. Then,
	\begin{align*}
	f^{G}_{v,w}=\tau_G\mathrm{dist}_\cT(i,j)+\tau_G\sum_{(x,y)\in P_{v,w}}\frac{f_{x,y}^{G_x}}{\tau_{G_x}},
	\end{align*}
	and so, 
	\begin{align}\label{2;temp1}
	r_{v,w}^G=\mathrm{dist}_\cT(i,j)+\sum_{(x,y)\in P_{v,w}}r_{x,y}^{G_x}.
	\end{align} 
	Moreover, 
	\begin{align}\label{2;temp2}
	\widehat{\mathbf{d}}_{G_i}^TR_{G}\be_{w}&=\mathbf{d}_{G_{i}}^TR_{G_{i}}\be_{w_1}+2\mathrm{dist}_\cT(i,j)m_{G_i}+2m_{G_i}\sum_{(x,y)\in P_{w_1,w}}r_{x,y}^{G_x}.
	\end{align}
\end{proposition}
\begin{proof}
	Let $l=\mathrm{dist}_\cT(i,j)+1$. In this proof, we shall use the same notation for the subgraphs and vertices used to set up $P_{v,w}$ with Figure \ref{Figure:P{v,w}}, so $P_{v,w}=\{(v_i,w_i)\,|\,i=1,\dots,l\}$ where $v=v_1$, and $w=w_l$. For $m=1,\dots,l$, we define $X_m$ as the set of $2$-tree spanning forests in $\mathcal{F}_{G}(v_1;w_l)$ such that the subtree with the vertex $v_1$ contains $v_m$ and the other with $w_l$ contains $w_m$. Similarly, for $m=1,\dots,l-1$, $Y_m$ is defined as the set of $2$-tree spanning forests in $\mathcal{F}_{G}(v_1;w_l)$ such that the subtree with $v_1$ contains $w_m$ and the other with $w_l$ contains $v_{m+1}$. Then, $X_m$'s and $Y_m$'s are mutually independent and $\mathcal{F}_{G}(v_1;w_l)=\left(\bigcup_{m=1}^l X_m\right)\cup \left(\bigcup_{m=1}^{l-1} Y_m\right)$. Applying an analogous argument as done in Proposition \ref{Proposition:dRd bridge}, we can find that $|X_m|=f_{v_m,w_m}^{G_{i_m}}\frac{\tau_G}{\tau_{G_{i_m}}}$ (note $G_{i_m}=G_{v_m}=G_{w_m}$) and $|Y_m|=\tau_G$. Hence, our desired results for $f^{G}_{v,w}$ and $r_{v,w}^G$ are established.
	
	We note that $P_{w_1,w}=\left(P_{v,w}\backslash\{(v,w_1)\}\right)\cup \{(w_1,w_1)\}$ and $r_{w_1,w_1}^{G_i}=0$. Using \eqref{2;temp1}, we have
	\begin{align*}
	\widehat{\mathbf{d}}_{G_{i}}^TR_{G}\be_{w_l}&=\sum\limits_{v\in V\left(G_i\right)}\mathrm{deg}_{G_i}(v)r_{v,w_l}^{G}\\
	&=\sum\limits_{v\in V\left(G_{i}\right)}\mathrm{deg}_{G_{i}}(v)\left(\mathrm{dist}_\cT(i,j)+r_{v,w_1}^{G_{i}}+\sum_{(x,y)\in P_{w_1,w}}r_{x,y}^{G_x}\right)\\
	&=\mathbf{d}_{G_{i}}^TR_{G_{i}}\be_{w_1}+2\mathrm{dist}_\cT(i,j)m_{G_{i}}+2m_{G_{i}}\sum_{(x,y)\in P_{w_1,w}}r_{x,y}^{G_x}.
	\end{align*}
\end{proof}

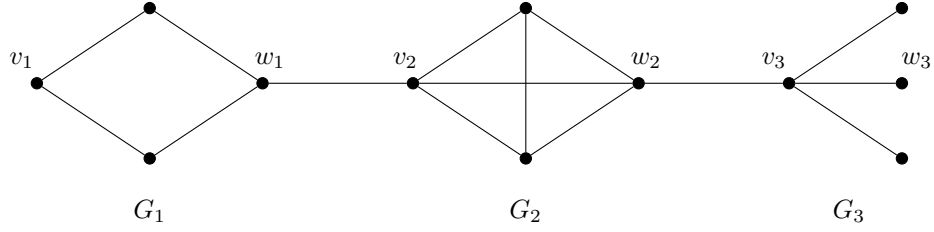
\begin{figure}[h!]
	\centering
	\begin{tikzpicture}
	\tikzset{enclosed/.style={draw, circle, inner sep=0pt, minimum size=.15cm, fill=black}}
	
	\node[enclosed, label={above, xshift=-0.2cm: $v_1$}] (v_1) at (-6.5,0) {};
	\node[enclosed, label={below, yshift=0cm: }] (v_2) at (-5,-1) {};
	\node[enclosed, label={above, yshift=0cm: }] (v_3) at (-5,1) {};
	\node[enclosed, label={above, xshift=0.1cm: $w_1$}] (v_4) at (-3.5,0) {};
	
	\node[enclosed, label={above, xshift=-0.1cm: $v_2$}] (v_5) at (-1.5,0) {};
	\node[enclosed, label={above, yshift=0cm: }] (v_6) at (0,1) {};
	\node[enclosed, label={right, xshift=0cm: }] (v_7) at (0,-1) {};
	\node[enclosed, label={above, xshift=0.1cm: $w_2$}] (v_8) at (1.5,0) {};
	
	\node[enclosed, label={above, xshift=-0.2cm: $v_3$}] (v_9) at (3.5,0) {};
	\node[enclosed, label={above, yshift=0cm: }] (v_10) at (5,1) {};
	\node[enclosed, label={above, xshift=0.2cm: $w_3$}] (v_11) at (5,0) {};
	\node[enclosed, label={right, xshift=0cm: }] (v_12) at (5,-1) {};
	
	\draw (v_1) -- (v_2); \draw (v_1) -- (v_3); 
	\draw (v_2) -- (v_4); \draw (v_3) -- (v_4);
	
	\draw (v_4) -- (v_5);
	
	\draw (v_5) -- (v_6); \draw (v_5) -- (v_7); \draw (v_5) -- (v_8);
	\draw (v_6) -- (v_7); \draw (v_6) -- (v_8); \draw (v_7) -- (v_8);
	
	\draw (v_8) -- (v_9);
	
	\draw (v_9) -- (v_10); \draw (v_9) -- (v_11); \draw (v_9) -- (v_12);
	
	\node[] at (-5,-1.7) {$G_1$};
	\node[] at (0,-1.7) {$G_2$};
	\node[] at (4.3,-1.7) {$G_3$};
	\end{tikzpicture}
	\caption{A chain of $3$ connected graphs with respect to a path on $3$ vertices, which is used in Examples \ref{Example:Cycle-Com-Star}, \ref{Example:Cycle-Com-Star2}, and \ref{Example:Cycle-Com-Star3}.}\label{Figure:cycle-com-star}
\end{figure}

We provide a straightforward example of these computations for the sake of clarity.
\begin{example}\label{Example:Cycle-Com-Star}
	Let $G_1$ be a cycle of length $4$, $G_2$ be a complete graph on $4$ vertices, and $G_3$ be a star on $4$ vertices. Suppose that $G$ is the chain of $G_1$, $G_2$, and $G_3$ with respect to a path on $3$ vertices as shown in Figure \ref{Figure:cycle-com-star}. Then, $\tau_{G_1}=4$, $\tau_{G_2}=16$ (by Cayley's formula (see \cite{book:GraphAndDigraphs})), $\tau_{G_3}=1$, and $\tau_G=64$. Moreover, one can find that we have
	\begin{align}\label{matrices:cycle-com-star}
	R_{G_1}=\frac{1}{4}\begin{bmatrix}
	0 & 3 & 4 & 3\\ 
	3 & 0 & 3 & 4\\ 
	4 & 3 & 0 & 3\\ 
	3 & 4 & 3 & 0
	\end{bmatrix},\;\;R_{G_2}=\frac{1}{16}\begin{bmatrix}
	0 & 8 & 8 & 8\\ 
	8 & 0 & 8 & 8\\ 
	8 & 8 & 0 & 8\\ 
	8 & 8 & 8 & 0
	\end{bmatrix},\;\;R_{G_3}=\frac{1}{1}\begin{bmatrix}
	0 & 1 & 1 & 1\\ 
	1 & 0 & 2 & 2\\ 
	1 & 2 & 0 & 2\\ 
	1 & 2 & 2 & 0
	\end{bmatrix},
	\end{align}
	where the first row and column in $R_{G_3}$ correspond to the vertex whose degree is $3$ in $G_3$. We have $f_{v_1,w_1}^{G_1}=4$, $f_{v_2,w_2}^{G_2}=8$. By Proposition \ref{Proposition:entries of f in chain}, $f_{v_1,w_3}^G=64(3-1)+64(\frac{4}{4}+\frac{8}{16}+1)=288$ and $r_{v_1,w_3}^G=\frac{1}{\tau_G}f_{v_1,w_3}^G=\frac{9}{2}$. Furthermore,
	\begin{align*}
	\widehat{\mathbf{d}}_{G_1}^TR_{G}\be_{w_3}	&=\mathbf{d}_{G_1}^TR_{G_1}\be_{w_1}+2(3-1)m_{G_1}+2m_{G_1}(r_{v_2,w_2}^{G_2}+r_{v_3,w_3}^{G_3})=33.
	\end{align*}
\end{example}

In order to compute $\kappa(G)$ when $G$ is a chain of connected graphs $G_1,\dots,G_k$ with respect to a tree $\cT$, we need to compute the quantity $\mathbf{d}_{G}^TR_{G}\mathbf{d}_G$. This quantity can be expressed as a weighted sum of the quantities $\mathbf{d}_{G}^TR_{G}\be_v$ for all $v\in V(G)$, and it is this term $\mathbf{d}_{G}^TR_{G}\be_v$ which we calculate in the next lemma, Lemma \ref{Lemma2:dRe}, using the results of Proposition \ref{Proposition:entries of f in chain}. As observed previously, since the graph $G$ has bridges, the resistance distances may be reduced and expressed in terms of resistance distances between vertices in the components $G_i$, with the expression taking into account the structure of the tree $\cT$. For any vertex $v\in V(G)$, there is a particular set of vertices in $G$ incident with the bridges of $G$ (corresponding to the edges of $\cT$) on which the reduction of our expression for $\mathbf{d}_{G}^TR_{G}\be_v$ will be based. This next observation makes the choice of these vertices clear, and its role is seen clearly in the statement of \Cref{Lemma2:dRe} and the proof of \Cref{Theorem:MainResult}.

\begin{observation}\label{Observation1}
	Let $\cT$ be a tree on $k\geq 2$ vertices and let $G$ be a chain of connected graphs $G_1,\dots,G_k$ with respect to $\cT$. For any pair of vertices $i$ and $j$ of $\cT$, they are disconnected upon the removal of any edge on the path from $i$ to $j$ in $\cT$. Hence, the deletion of one of the bridges in $\cB_G$ that correspond to the edges on the path from $i$ and $j$ in $\cT$ results in two components, where one contains $G_i$ and the other contains $G_j$. In particular, among those bridges there exists exactly one bridge such that it is incident to some vertex $w$ of $G_i$; in other words, $w$ is the only vertex of $G_i$ such that it is incident to a bridge in $\cB_G$ and the deletion of that bridge results in two components, where one contains $G_{i}$ and the other contains $G_{j}$.
	
	Given a vertex $v\in V(G)$, there exists $i_0\in\{1,\dots,k\}$ such that $v\in G_{i_0}$. From the above, there are exactly $(k-1)$ vertices $w_l$ for $l\in\{1,\dots,k\}\backslash\{i_0\}$ such that $w_l\in V(G_l)$, and $w_l$ is incident to a bridge in $\cB_G$ whose deletion results in two components, where the one with $w_l$ contains $G_l$ and the other does $G_{i_0}$. Furthermore, for two distinct vertices $w_{l_1}$ and $w_{l_2}$, the corresponding bridges are distinct---that is, the $(k-1)$ vertices correspond to the bridges in $\cB_G$, and vice versa. Fixing any vertex $v\in V(G)$ uniquely determines a set of vertices which are in one-to-one correspondence with the bridges in $\cB_G$.
\end{observation}

Here we introduce some notation. Let $G$ be a chain of connected graphs $G_1,\dots,G_k$ with respect to a tree $\cT$ on $k$ vertices. Let $\widetilde{V}_k^G$ denote the set of the vertices in $V(G_k)$ that are incident to bridges in $\cB_G$. Note that $|\widetilde{V}_k^G|=\mathrm{deg}_\cT(k)$. For each bridge $x\sim y\in\cB_G$, we have exactly two components in $G\backslash x\sim y$. We use $W_x^G$ (resp. $W_y^G$) to denote the number of edges of the component with $x$ (resp. $y$) in $G\backslash x\sim y$. Let $\widebar{W}_x^G:=m_G-W_x^G$ and $\widebar{W}_y^G:=m_G-W_y^G$. Note that $\widebar{W}_x^G=W_y^G+1$ and $\widebar{W}_y^G=W_x^G+1$. We shall omit the superscript $G$ if it is clear from the context. In this article, the superscript $G$ for $\widetilde{V}_k^G$, $W_x^G$ and $\widebar{W}_x^G$ is only used in the proof of Theorem \ref{Theorem:MainResult} and in Question \ref{min-max question1}.

In the next lemma, we give a formula for $\mathbf{d}_G^TR_{G}\be_{w_k}$ for any vertex $w_k \in V(G)$ using the results of \Cref{Proposition:entries of f in chain} and the notation outlined so far.
\begin{lemma}\label{Lemma2:dRe}
	Given a tree $\cT$ on $k$ vertices and connected graphs $G_1,\dots,G_k$, let $G$ be a chain of $G_1,\dots,G_k$ with respect to $\cT$. Let $w_k\in V(G_k)$. From Observation \ref{Observation1}, for $i=1,\dots,k-1$, we may define $w_i$ to be the unique vertex of $G_i$ such that $w_i$ is incident to a bridge in $\cB_G$ whose deletion results in two components, where the one with $w_i$ contains $G_i$ and the other contains $G_k$. Then,
	\begin{align}\label{2;temp0}
	\mathbf{d}_G^TR_{G}\be_{w_k}=\sum_{i=1}^{k}\mathbf{d}_{G_i}^TR_{G_i}\be_{w_i}+\sum_{i=1}^{k}\sum_{z\in \widetilde{V}_i}2\widebar{W}_zr_{z,w_i}^{G_i}+\sum_{i=1}^{k-1}(2W_{w_i}+1).
	\end{align}
\end{lemma}
\begin{proof}
	We note that $$\mathbf{d}_G=\sum_{j=1}^k\widehat{\mathbf{d}}_{G_j}+\sum_{x\sim y\in \cB_G}(\be_x+\be_y)=\sum_{j=1}^k\widehat{\mathbf{d}}_{G_j}+\sum_{j=1}^k\sum_{z\in \widetilde{V}_j}\be_z.$$
	Using \eqref{2;temp1} and \eqref{2;temp2}, we have
	\begin{align}\nonumber
	\mathbf{d}_G^TR_{G}\be_{w_k}=&\sum_{i=1}^k\widehat{\mathbf{d}}_{G_i}^TR_{G}\be_{w_k}+\sum_{i=1}^k\sum_{z\in \widetilde{V}_i}r_{z,w_k}^{G}\\\nonumber
	=&\sum_{i=1}^k\left(\mathbf{d}_{G_{i}}^TR_{G_{i}}\be_{w_i}+2\mathrm{dist}_\cT(i,k)m_{G_{i}}+2m_{G_{i}}\sum_{(x,y)\in P_{w_i,w_k}}r_{x,y}^{G_x}\right)\\\nonumber
	&\quad+\sum_{i=1}^k\sum_{z\in \widetilde{V}_i}\left(\mathrm{dist}_\cT(i,k)+\sum_{(x,y)\in P_{z,w_k}}r_{x,y}^{G_x}\right)\\\label{2;temp00}
	=&\sum_{i=1}^k\left(2m_{G_{i}}\sum_{(x,y)\in P_{w_i,w_k}}r_{x,y}^{G_x}+\sum_{z\in \widetilde{V}_i}\sum_{(x,y)\in P_{z,w_k}}r_{x,y}^{G_x}\right)\\\label{2;temp01}
	&\quad+\sum_{i=1}^k\mathbf{d}_{G_{i}}^TR_{G_{i}}\be_{w_i}+\sum_{i=1}^k\mathrm{dist}_\cT(i,k)(2m_{G_i}+\mathrm{deg}_\cT(i)).
	\end{align}
	
	First, we claim that the right side in \eqref{2;temp00} can be written as the second summand in the right side of \eqref{2;temp0}. In \eqref{2;temp00}, we note that for any pair $(x,y)$ in $P_{w_i,w_k}\cup P_{z,w_k}$, we have $x,y\in \widetilde{V}_j$ for some $1\leq j\leq k$. Thus, in order to establish the claim, it suffices to show that given $1\leq s\leq k$, for each $\hat{z}\in\widetilde{V}_s$, the coefficient of $r_{\hat{z},w_s}^{G_s}$ in the expansion of the summation in \eqref{2;temp00} is $2\widebar{W}_{\hat{z}}$. Since $r_{w_s,w_s}^{G_s}=0$, we may assume $\hat{z}\neq w_s$. For the bridge $b$ in $\cB_G$ incident to $\hat{z}$, we let $C$ be the component in $G\backslash b$ such that $\hat{z}\notin V(C)$. Consider the expansion of $2m_{G_{i}} \sum_{(x,y)\in P_{w_i,w_k}}r_{x,y}^{G_x}$. We can find from Observation \ref{Observation2} that $(\hat{z},w_s)\in P_{w_i,w_k}$ so that $r_{\hat{z},w_s}^{G_s}$ appears with coefficient $2m_{G_i}$ if for $i\neq s$, vertex $s$ lies on the path from $i$ to $k$ in $\cT$, and $w_i\in V(C)$. Similarly, for the expansion of $\sum_{z\in \widetilde{V}_i}\sum_{(x,y)\in P_{z,w_k}} r_{x,y}^{G_x}$, $(\hat{z},w_s)\in P_{z,w_k}$ so that $r_{\hat{z},w_s}^{G_s}$ appears with coefficient $\mathrm{deg}_\cT(i)$ if for $i\neq s$, vertex $s$ lies on the path from $i$ to $k$ in $\cT$, and $z\in V(C)$; and $r_{\hat{z},w_s}^{G_s}$ appears with coefficient $1$ if $i=s$ and $z=\hat{z}$. Note that $\mathrm{deg}_\cT(i)$ is equal to the number of bridges in $\cB_G$ that are  incident to some vertices of $G_i$. It follows that $r_{\hat{z},w_s}^{G_s}$ has coefficient $2\widebar{W}_{\hat{z}}$. Hence, our claim is established.
	
	Now, we consider the last summand in $\eqref{2;temp01}$; and since $\mathrm{dist}_{\cT}(k,k)=0$, we shall show that
	\begin{align*}
	\sum_{i=1}^{k-1}(2W_{w_i}+1)=\sum_{i=1}^{k-1} \mathrm{dist}_\cT(i,k)(2m_{G_i}+\mathrm{deg}_\cT(i)).
	\end{align*}
	For $i=1,\dots,k-1$, $W_{w_i}$ is the number of edges of the component $H$ with $w_i$ in $G\backslash e$, where $e$ is the bridge in $\cB_G$ incident to $w_i$. Then, $w_k\notin V(H)$. We note that $H$ consists of some graphs $G_{j_1},\dots,G_{j_l}$ for some $l\geq 1$ and some bridges in $\cB_G$ so that the number of bridges in $\cB_G$ that belong to $H$ is $l-1$. We further note that for $r=1,\dots,l$, the number of bridges in $\cB_G$ incident to some vertices of $G_{j_r}$ is $\mathrm{deg}_\cT(j_r)$. Letting $j_l=i$, we have $2(l-1)=\sum_{r=1}^{l-1}\mathrm{deg}_\cT(j_r)+(\mathrm{deg}_\cT(i)-1)$. Hence, $2W_{w_i}+1$ can be written in terms of $m_{G_{j_1}},\dots,m_{G_{j_l}}$ and $\mathrm{deg}_\cT(j_1),\cdots,\mathrm{deg}_\cT(j_l)$. Therefore, it is enough to show that in computation of $\sum_{i=1}^{k-1}(2W_{w_i}+1)$, for each $s=1,\dots, k-1$, $m_{G_s}$ and $\mathrm{deg}_\cT (s)$ appear exactly $2\mathrm{dist}_\cT(s,k)$ and $\mathrm{dist}_\cT(s,k)$ times, respectively.
	
	Note that $i\neq k$. If vertex $i$ lies on the path from $s$ to $k$ in $\cT$, then $V(G_s)\subseteq V(H)$ and so $m_{G_s}$ appears exactly once in computation of $W_{w_i}$. This implies that $m_{G_s}$ appears $2\mathrm{dist}_\cT(s,k)$ times in computation of $\sum_{i=1}^{k-1}(2W_{w_i}+1)$. 
	
	If for $i\neq s$, vertex $i$ lies on the path from $s$ to $k$, then $\mathrm{deg}_\cT(s)$ appears exactly once in computation of $2W_{w_i}$; and, if $i=s$, then $\mathrm{deg}_\cT(s)-1$ appears once in computation of $2W_{w_s}$. It follows that $\mathrm{deg}_\cT(s)$ appears $\mathrm{dist}_\cT(s,k)$ times in computation of $\sum_{i=1}^{k-1}(2W_{w_i}+1)$.
\end{proof}

\begin{example}\label{Example:Cycle-Com-Star2}
	Maintaining Example \ref{Example:Cycle-Com-Star} with the same notation, let us find $\mathbf{d}_G^TR_G\be_{w_3}$ with Lemma \ref{Lemma2:dRe}. Note $\cB_G=\{w_1\sim v_2, w_2\sim v_3\}$, $v_1\notin \widetilde{V}_1$, and $w_3\notin \widetilde{V}_3$. We have $\widebar{W}_{v_2}=5$, $\widebar{W}_{v_3}=12$, $W_{w_1}=4$, and $W_{w_2}=11$. Then,
	\begin{align*}
	\mathbf{d}_G^TR_G\be_{w_3}=&\sum_{i=1}^{3}\mathbf{d}_{G_i}^TR_{G_i}\be_{w_i}+2\widebar{W}_{v_2}r_{v_2,w_2}^{G_2}+2\widebar{W}_{v_3}r_{v_3,w_3}^{G_3}\\
	&+(2W_{w_1}+1)+(2W_{w_2}+1)\\
	=&(5+4.5+7)+5+24+9+23=77.5.
	\end{align*}
\end{example}

Here is our main result of this article.

\begin{theorem}\label{Theorem:MainResult}
	Let $\cT$ be a tree on $k$ vertices, and let $G$ be a chain of connected graphs $G_1,\dots,G_k$ with respect to $\cT$. Then, 
	\begin{align}\label{2;formula K(G)}
	\begin{split}
	\kappa(G)=&\sum\limits_{i=1}^{k}\frac{m_{G_i}}{m_G}\kappa(G_i)+\hspace*{0cm}\sum\limits_{x\sim y \in \cB_G}\left(\frac{\widebar{W}_x}{m_G}\mathbf{d}_{G_x}^TR_{G_x}\be_x+\frac{\widebar{W}_y}{m_G}\mathbf{d}_{G_y}^TR_{G_y}\be_y\right)\\
	&+\sum\limits_{i=1}^{k}\sum\limits_{\small(z_1,z_2) \in \widetilde{V}_i\times\widetilde{V}_i\normalsize}\hspace*{-0.5cm}\frac{\widebar{W}_{z_1}\widebar{W}_{z_2}}{m_G}r_{z_1,z_2}^{G_i}\hspace*{-0.05cm}+\hspace*{-0.25cm}\sum\limits_{x\sim y \in \cB_G}\hspace*{-0.3cm}\frac{(2\widebar{W}_x-1)(2\widebar{W}_y-1)}{2m_G}.
	\end{split}
	\end{align}
\end{theorem}
\begin{proof}
	We shall use induction on $k$ for this proof. Clearly, the statement holds for $k=1$. Suppose that for $k\geq 1$, this expression holds for every chain of connected graphs with respect to any tree on $k$ vertices. Let $\widetilde{\cT}$ be a tree with vertex set $\{1,\dots,k+1\}$. Without loss of generality, assume $\mathrm{deg}_{\widetilde{\cT}}(k+1)=1$ and $k\sim (k+1)$ is an edge of $\widetilde{\cT}$. Let $\widetilde{G}$ be a chain of connected graphs $G_1,\dots,G_{k+1}$ with respect to $\widetilde{\cT}$. We first state what we need to obtain to complete our induction:
	\begin{align}\label{what we need to obtain}
	\begin{split}
	\kappa(\widetilde{G})=&\sum\limits_{i=1}^{k+1}\frac{m_{G_i}}{m_{\widetilde{G}}}\kappa(G_i)+
	\hspace*{0cm}\sum\limits_{x\sim y \in \cB_{\widetilde{G}}}\left(\frac{\widebar{W}_x^{\widetilde{G}}}{m_{\widetilde{G}}}\mathbf{d}_{G_x}^TR_{G_x}\be_x+\frac{\widebar{W}_y^{\widetilde{G}}}{m_{\widetilde{G}}}\mathbf{d}_{G_y}^TR_{G_y}\be_y\right)\\
	&+\sum\limits_{i=1}^{k+1}\sum\limits_{\small(z_1,z_2) \in \widetilde{V}_i^{\widetilde{G}}\times\widetilde{V}_i^{\widetilde{G}}\normalsize}\hspace*{-0.5cm}\frac{\widebar{W}_{z_1}^{\widetilde{G}}\widebar{W}_{z_2}^{\widetilde{G}}}{m_{\widetilde{G}}}r_{z_1,z_2}^{G_i}\hspace*{-0.05cm}
	+\hspace*{-0.25cm}\sum\limits_{x\sim y \in \cB_{\widetilde{G}}}\hspace*{-0.3cm}\frac{(2\widebar{W}_x^{\widetilde{G}}-1)(2\widebar{W}_y^{\widetilde{G}}-1)}{2m_{\widetilde{G}}}.
	\end{split}
	\end{align}
	
	Now we shall derive \eqref{what we need to obtain} from the inductive hypothesis and the results we have obtained. Suppose that $G$ is the graph obtained from $\widetilde{G}$ by removing all vertices of $G_{k+1}$ and all edges incident to vertices of $G_{k+1}$. Further, assume that $w\sim v$ is the bridge in $\cB_{\widetilde{G}}$ where $w\in V(G_k)$ and $v\in V(G_{k+1})$. Then, $G$ is a chain of $G_1,\dots,G_k$ with respect to $\cT$, where $\cT$ is the tree obtained from $\widetilde{\cT}$ by removing the pendent vertex $k+1$ and the edge $k\sim (k+1)$. Using \Cref{Prop1:bridge between 2}, we have
	\begin{align}\label{2;thm temp 1}
	\begin{split}
	\kappa(\widetilde{G})=&\frac{m_{G}}{m_{\widetilde{G}}}\kappa(G)+\frac{m_{G_{k+1}}}{m_{\widetilde{G}}}\kappa(G_{k+1})+\frac{m_{G_{k+1}}+1}{m_{\widetilde{G}}}\mathbf{d}_{G}^TR_{G}\be_{w}\\
	&+\frac{m_G+1}{m_{\widetilde{G}}}\mathbf{d}_{G_{k+1}}^TR_{G_{k+1}}\be_{v}+\frac{(2m_G+1)(2m_{G_{k+1}}+1)}{2m_{\widetilde{G}}}.
	\end{split}
	\end{align}
	
	Let $w_k=w$. By \Cref{Observation1}, for $i=1,\dots,k-1$, we may define $w_i$ to be the vertex of $G_i$ such that $w_i$ is incident to a bridge in $\cB_G$ whose deletion results in two components, where the one with $w_i$ contains $G_i$ and the other contains $G_k$; in this sense, there is one-to-one correspondence between $\cB_G$ and the set of $w_1,\dots,w_{k-1}$. Here we remark the following, which is used in several places of this proof.
	\begin{enumerate}[label=(R)]
		\item\label{temporary remark} In the same sense above, there is one-to-one correspondence between $\cB_{\widetilde{G}}$ and the set of $w_1,\dots,w_{k}$. Note that $k$ is adjacent to $k+1$ in $\widetilde{\cT}$. For $i=1,\dots,k$, after the deletion of the bridge corresponding to $w_i$ from $\widetilde{G}$, the component not having $w_i$ contains all vertices of $G_{k+1}$.
	\end{enumerate}
	
	In \eqref{2;thm temp 1}, applying the inductive hypothesis to $\kappa(G)$ and \Cref{Lemma2:dRe} to $\mathbf{d}_{G}^TR_{G}\be_{w}$, we obtain     
	\begin{align}\label{2;thm temp 2}
	\begin{split}
	\frac{m_{G}}{m_{\widetilde{G}}}\kappa(G)=&\sum\limits_{i=1}^{k}\frac{m_{G_i}}{m_{\widetilde{G}}}\kappa(G_i)+\sum\limits_{x\sim y \in \cB_G}\left(\frac{\widebar{W}_x^G}{m_{\widetilde{G}}}\mathbf{d}_{G_x}^TR_{G_x}\be_x+\frac{\widebar{W}_y^G}{m_{\widetilde{G}}}\mathbf{d}_{G_y}^TR_{G_y}\be_y\right)\\
	&+\sum\limits_{i=1}^{k}\sum\limits_{\small(z_1,z_2) \in \widetilde{V}_i^G\times\widetilde{V}_i^G\normalsize}\hspace*{-0.3cm}\frac{\widebar{W}_{z_1}^G\widebar{W}_{z_2}^G}{m_{\widetilde{G}}}r_{z_1,z_2}^{G_i}+\hspace*{-0.25cm}\sum\limits_{x\sim y \in \cB_G}\hspace*{-0.25cm}\frac{(2\widebar{W}_x^G-1)(2\widebar{W}_y^G-1)}{2m_{\widetilde{G}}},
	\end{split}
	\end{align}
	and
	\begin{align}\label{2;thm temp 3}
	\begin{split}
	\frac{m_{G_{k+1}}\hspace*{-0.1cm}+1}{m_{\widetilde{G}}}\mathbf{d}_{G}^TR_{G}\be_{w}=&\frac{m_{G_{k+1}}\hspace*{-0.1cm}+1}{m_{\widetilde{G}}}\left(\sum_{i=1}^{k}\mathbf{d}_{G_i}^TR_{G_i}\be_{w_i}\hspace*{-0.1cm}+\sum_{i=1}^{k}\sum_{z\in \widetilde{V}_i^G}2\widebar{W}_z^Gr_{z,w_i}^{G_i}\right)\\
	&+\frac{m_{G_{k+1}}\hspace*{-0.1cm}+1}{m_{\widetilde{G}}}\sum_{i=1}^{k-1}(2W_{w_i}^G+1).
	\end{split}
	\end{align}	
	
	Now, we consider \eqref{2;thm temp 1} together with \eqref{2;thm temp 2} and \eqref{2;thm temp 3} in order to derive \eqref{what we need to obtain}. First, $\sum\limits_{i=1}^{k+1}\frac{m_{G_i}}{m_{\widetilde{G}}}\kappa(G_i)$ can be obtained from the sum of $\frac{m_{G_{k+1}}}{m_{\widetilde{G}}}\kappa(G_{k+1})$ in \eqref{2;thm temp 1} and $\sum\limits_{i=1}^{k}\frac{m_{G_i}}{m_{\widetilde{G}}}\kappa(G_i)$ in \eqref{2;thm temp 2}. Next, we claim that the sum of $\frac{m_{G_{k+1}}+1}{m_{\widetilde{G}}}\sum_{i=1}^{k}\mathbf{d}_{G_i}^TR_{G_i}\be_{w_i}$ in \eqref{2;thm temp 3}, $\sum\limits_{x\sim y \in \cB_G}\left(\frac{\widebar{W}_x^G}{m_{\widetilde{G}}}\mathbf{d}_{G_x}^TR_{G_x}\be_x+\frac{\widebar{W}_y^G}{m_{\widetilde{G}}}\mathbf{d}_{G_y}^TR_{G_y}\be_y\right)$ in \eqref{2;thm temp 2}, and $\frac{m_G+1}{m_{\widetilde{G}}}\mathbf{d}_{G_{k+1}}^TR_{G_{k+1}}\be_{v}$ in \eqref{2;thm temp 1} yields 
	$$
	\sum\limits_{x\sim y \in \cB_{\widetilde{G}}}\left(\frac{\widebar{W}_x^{\widetilde{G}}}{m_{\widetilde{G}}}\mathbf{d}_{G_x}^TR_{G_x}\be_x+\frac{\widebar{W}_y^{\widetilde{G}}}{m_{\widetilde{G}}}\mathbf{d}_{G_y}^TR_{G_y}\be_y\right).
	$$ 
	In order to establish the claim, we consider two cases: either a bridge $x\sim y$ belongs to $\cB_G$ or it is $w\sim v$. First, let us choose $x\sim y \in \cB_G$. Note that $|\cB_G|=k-1$. By \ref{temporary remark}, one of two components in $\widetilde{G}\backslash x\sim y$ must contain all vertices of $G_{k+1}$, say the component containing $y$ does so. Then, $x=w_{i_0}$ for some $i_0\in\{1,\dots k-1\}$. Moreover, $\widebar{W}_{y}^{\widetilde{G}}=\widebar{W}_{y}^G$ and $\widebar{W}_{x}^{\widetilde{G}}=\widebar{W}_{x}^G+m_{G_{k+1}}+1$. Hence, we obtain
	\begin{align*}
	\frac{\widebar{W}_x^{\widetilde{G}}}{m_{\widetilde{G}}}\mathbf{d}_{G_x}^TR_{G_x}\be_x+\frac{\widebar{W}_y^{\widetilde{G}}}{m_{\widetilde{G}}}\mathbf{d}_{G_y}^TR_{G_y}\be_y
	& =\frac{\widebar{W}_x^G}{m_{\widetilde{G}}}\mathbf{d}_{G_x}^TR_{G_x}\be_x+\frac{\widebar{W}_y^G}{m_{\widetilde{G}}}\mathbf{d}_{G_y}^TR_{G_y}\be_y\\
	&\qquad\quad +\frac{m_{G_{k+1}}+1}{m_{\widetilde{G}}}\mathbf{d}_{G_{i_0}}^TR_{G_{i_0}}\be_{w_{i_0}}.
	\end{align*}
	Consider the latter case $w\sim v$ (note $w=w_k$). Considering \ref{temporary remark}, we have $$\hspace*{-0.05cm}\frac{\widebar{W}_w^{\widetilde{G}}}{m_{\widetilde{G}}}\mathbf{d}_{G_w}^TR_{G_w}\be_w+\frac{\widebar{W}_v^{\widetilde{G}}}{m_{\widetilde{G}}}\mathbf{d}_{G_v}^TR_{G_v}\be_v=\frac{m_{G_{k+1}}\hspace*{-0.1cm}+1}{m_{\widetilde{G}}}\mathbf{d}_{G_{k}}^TR_{G_{k}}\be_{w_{k}}+\frac{m_G+1}{m_{\widetilde{G}}}\mathbf{d}_{G_{k+1}}^TR_{G_{k+1}}\be_{v}.$$
	We note again that the bridges in $\cB_{\widetilde{G}}$ correspond to $w_1,\dots,w_k$. Therefore, our claim is established, as desired.
	
	After that, we shall show that 
	\begin{align}\label{2;temp 4}
	\begin{split}
	\sum\limits_{i=1}^{k+1}\sum\limits_{\small(z_1,z_2) \in \widetilde{V}_i^{\widetilde{G}}\times\widetilde{V}_i^{\widetilde{G}}\normalsize}\hspace*{-0.4cm}\frac{\widebar{W}_{z_1}^{\widetilde{G}}\widebar{W}_{z_2}^{\widetilde{G}}}{m_{\widetilde{G}}}r_{z_1,z_2}^{G_i}
	& = \sum\limits_{i=1}^{k}\sum\limits_{\small(z_1,z_2) \in \widetilde{V}_i^G\times\widetilde{V}_i^G\normalsize}\hspace*{-0.3cm}\frac{\widebar{W}_{z_1}^G\widebar{W}_{z_2}^G}{m_{\widetilde{G}}}r_{z_1,z_2}^{G_i}\\
	& \qquad \qquad +\frac{m_{G_{k+1}}\hspace*{-0.1cm}+1}{m_{\widetilde{G}}}\sum_{i=1}^{k}\sum_{z\in \widetilde{V}_i^G}2\widebar{W}_z^Gr_{z,w_i}^{G_i}.
	\end{split}
	\end{align}
	In order to do that, we shall compare coefficients for both sides. We note $r_{a,a}^{H}=0$ for any graph $H$ and vertex $a$ of $H$. Since $\mathrm{deg}_{\widetilde{\cT}}(k+1)=1$, $\widetilde{V}_{k+1}^{\widetilde{G}}$ is a singleton. So, it is enough to show that given $i=1,\dots,k$, for any $(z_1,z_2)\in \widetilde{V}_i^{\widetilde{G}}\times\widetilde{V}_i^{\widetilde{G}}$ with $z_1\neq z_2$, the respective coefficients of $r_{z_1,z_2}^{G_i}$ on both sides in \eqref{2;temp 4} are equal. 
	
	For $(z_1,z_2)\in \widetilde{V}_i^{\widetilde{G}}\times\widetilde{V}_i^{\widetilde{G}}$ with $z_1\neq z_2$, we let $e_1$ and $e_2$ be the bridges in $\cB_{\widetilde{G}}$ such that they are incident to $z_1$ and $z_2$, respectively. Note that $w_i\in \widetilde{V}_i^{\widetilde{G}}$. We now consider two cases: \begin{enumerate*}[label=(\roman*)]
		\item\label{2;case11} one of $z_1$ and $z_2$ is $w_i$, and 
		\item\label{2;case12} neither of them is $w_i$.
	\end{enumerate*}
	For case \ref{2;case11}, we suppose without loss of generality that $z_2=w_i$. From \ref{temporary remark}, the component without $w_i$ in $\widetilde{G}\backslash e_2$ contains $G_{k+1}$, while the component without $z_1$ in $\widetilde{G}\backslash e_1$ does not contain  $G_{k+1}$. Hence, $\widebar{W}_{z_1}^{\widetilde{G}}=\widebar{W}_{z_1}^G$ and $\widebar{W}_{w_i}^{\widetilde{G}}=\widebar{W}_{w_i}^G+m_{G_{k+1}}+1$ for $i=1,\dots,k-1$; and $\widebar{W}_{z_1}^{\widetilde{G}}=\widebar{W}_{z_1}^G$ and $\widebar{W}_{w_k}^{\widetilde{G}}=m_{G_{k+1}}+1$. We note that $(z_1,w_i), (w_i,z_1)\in \widetilde{V}_i^{\widetilde{G}}\times\widetilde{V}_i^{\widetilde{G}}$ and $r_{z_1,w_i}^{G_i}=r_{w_i,z_1}^{G_i}$. For $i=1,\dots,k-1$, considering that $\widetilde{V}_i^G\times\widetilde{V}_i^G$ contains two distinct elements $(z_1,w_i)$ and $(w_i,z_1)$, replacing $\widebar{W}_{z_1}^{\widetilde{G}}$ and $\widebar{W}_{z_2}^{\widetilde{G}}$ in the left side of \eqref{2;temp 4} by $\widebar{W}_{z_1}^G$ and $\widebar{W}_{w_i}^G+m_{G_{k+1}}+1$, respectively, one can find that the respective coefficients of $r_{z_1,w_i}^{G_i}$ on both sides in \eqref{2;temp 4} are equal. Similarly, for the case $i=k$, noting $(z_1,w_i),(w_i,z_1)\notin \widetilde{V}_i^G\times\widetilde{V}_i^G$, the equality for the respective coefficients of $r_{z_1,w_i}^{G_i}$ can be obtained. Let us consider \ref{2;case12}. Then, for $j=1,2$, the component without $z_j$ in $\widetilde{G}\backslash e_j$ does not contain  $G_{k+1}$. So, $\widebar{W}_{z_1}^{\widetilde{G}}=\widebar{W}_{z_1}^G$ and $\widebar{W}_{z_2}^{\widetilde{G}}=\widebar{W}_{z_2}^G$. Hence, the respective coefficients of $r_{z_1,z_2}^{G_i}$ on both sides in \eqref{2;temp 4} are equal. Therefore, the equality in \eqref{2;temp 4} follows.
	
	Finally, for the completion of the proof, we only need to show
	\begin{align}\label{2;temp5}
	\begin{split}
	&\sum\limits_{x\sim y \in \cB_{\widetilde{G}}}\frac{(2\widebar{W}_x^{\widetilde{G}}-1)(2\widebar{W}_y^{\widetilde{G}}-1)}{2m_{\widetilde{G}}}\\
	=&\sum\limits_{x\sim y \in \cB_G}\hspace*{-0.25cm}\frac{(2\widebar{W}_x^G-1)(2\widebar{W}_y^G-1)}{2m_{\widetilde{G}}}
	+\frac{(2m_G+1)(2m_{G_{k+1}}+1)}{2m_{\widetilde{G}}}\\
	&+\frac{2(m_{G_{k+1}}\hspace*{-0.1cm}+1)}{2m_{\widetilde{G}}}\sum_{i=1}^{k-1}(2W_{w_i}+1).
	\end{split}
	\end{align}
	Recall \ref{temporary remark}. Choose a bridge $x\sim y \in \cB_{\widetilde{G}}$. If $x\sim y$ is $w\sim v$, then $\widebar{W}_w^{\widetilde{G}}=m_{G_{k+1}}+1$ and $\widebar{W}_v^{\widetilde{G}}=m_{G}+1$. Suppose that $x\sim y$ is in $\cB_G$. One of two components in $\widetilde{G}\backslash x\sim y$ must contain $G_{k+1}$, say the component with $y$ does so. Then, $\widebar{W}_x^{\widetilde{G}}=\widebar{W}_x^G+m_{G_{k+1}}+1$ and $\widebar{W}_y^{\widetilde{G}}=\widebar{W}_y^G$. Moreover, $x$ must be $w_{i_0}$ for some $1\leq i_0\leq k-1$, and so $W_{w_{i_0}}=\widebar{W}_y^G-1$. Then, one can verify that substituting the above expressions appropriately for $\widebar{W}_x^{\widetilde{G}}$ and $\widebar{W}_y^{\widetilde{G}}$ for $x\sim y \in \cB_{\widetilde{G}}$ in the left side of \eqref{2;temp5} yields the right side.  
\end{proof}

\begin{remark}\label{Remakr:V_i for pendent}
	Continuing the notation and result in Theorem \ref{Theorem:MainResult}, if $\widetilde{V}_{i_0}$ contains exactly one element for some $1\leq i_0\leq k$, \textit{i.e,} $\mathrm{deg}_\cT(i_0)=1$, then $r_{z,z}^{G_{i_0}}=0$ for $(z,z)\in \widetilde{V}_{i_0}\times \widetilde{V}_{i_0}$. Hence, when it comes to computation of $\sum\limits_{\small(z_1,z_2) \in \widetilde{V}_i\times\widetilde{V}_i\normalsize}\hspace*{-0.5cm}\frac{\widebar{W}_{z_1}\widebar{W}_{z_2}}{m_G}r_{z_1,z_2}^{G_i}$, we only need to consider indices $i$ such that $\mathrm{deg}_\cT(i)>1$.
\end{remark}

\begin{example}\label{Example:Cycle-Com-Star3}
	Continuing Examples \ref{Example:Cycle-Com-Star} and \ref{Example:Cycle-Com-Star2} with the same notation, we shall obtain $\kappa(G)$ through Theorem \ref{Theorem:MainResult}. It can be found that $\kappa(G_1)=2.5$, $\kappa(G_2)=2.25$, and $\kappa(G_3)=7.5$. Note that $\widetilde{V}_1$ and $\widetilde{V}_3$ both contain a single element. Then, we have
	\begin{align*}
	m_G\kappa(G)=&\sum\limits_{i=1}^{3}m_{G_i}\kappa(G_i)+\widebar{W}_{w_1}\mathbf{d}_{G_1}^TR_{G_1}\be_{w_1}+\widebar{W}_{v_2}\mathbf{d}_{G_2}^TR_{G_2}\be_{v_2}+\widebar{W}_{w_2}\mathbf{d}_{G_2}^TR_{G_2}\be_{w_2}\\
	&\quad +\widebar{W}_{v_3}\mathbf{d}_{G_3}^TR_{G_3}\be_{v_3}+\widebar{W}_{v_2}\widebar{W}_{w_2}r_{v_2,w_2}^{G_2}+\widebar{W}_{w_2}\widebar{W}_{v_2}r_{w_2,v_2}^{G_2}\\
	& \quad\quad +\frac{1}{2}(2\widebar{W}_{w_1}-1)(2\widebar{W}_{v_2}-1)+\frac{1}{2}(2\widebar{W}_{w_2}-1)(2\widebar{W}_{v_3}-1).
	\end{align*}
	One can verify that $\kappa(G)=\frac{357.5}{15}= 23.8\dot{3}$.
\end{example}

The formula for $\kappa(G)$ in \Cref{Theorem:MainResult} derives its importance from the fact that understanding the resistance matrices and the degree vectors of $G_1,\dots,G_k$ allows the computation of $\kappa(G)$ without calculating this quantities for $G$ from scratch. 

\begin{remark}\label{Remark:complexity}
	For computation of the Moore--Penrose inverse, the singular value decomposition is used; so, for an $m\times m$ matrix $A$, the cost for computation of $A^\dagger$ is $\mathcal{O}(m^3)$ \cite{hogben2013handbook} where $\mathcal{O}$ stands for the big O notation. Continuing Theorem \ref{Theorem:MainResult}, we suppose that $|V(G)|=n$ and $|V(G_i)|=n_i$ for $i=1,\dots,k$. Then, it follows that the cost for computation of the left side of \eqref{2;formula K(G)} is $\mathcal{O}(n^3)$; and the cost of the right side of \eqref{2;formula K(G)} is $\mathcal{O}(n_1^3+\cdots+n_k^3)$.
\end{remark}

\begin{remark}
	Given a chain $G$ of trees $\cT_1,\dots,\cT_k$ with respect to some tree on $k$ vertices, the formula for $\kappa(G)$ in \Cref{Theorem:MainResult} is equivalent to that in Proposition 2.2 of \cite{ciardo2020kemeny}.
\end{remark}

In the following example, we provide a formula for Kemeny's constant for a chain of connected graphs $G_1,\dots,G_k$ with respect to a star on $k$ vertices, as described in Figure \ref{Figure:Chain of Graphs}.
\begin{example}
	Let $k\geq 3$, and $\cT$ be a star with vertex set $\{1,\dots, k\}$. Suppose that vertex $1$ is of degree $k-1$. Let $G$ be a chain of connected graphs $G_1,\dots,G_{k}$ with respect to $\cT$. We may assume that $\cB_G=\{x_i\sim y_i| i=1,\dots,k-1\}$, $x_i\in V(G_1)$ and $y_i\in V(G_{i+1})$ for $1\leq i\leq k-1$. Then, we have $\widebar{W}_{x_i}=m_{G_{i+1}}+1$ and $\widebar{W}_{y_i}=m_G-m_{G_{i+1}}$. By Theorem \ref{Theorem:MainResult} with Remark \ref{Remakr:V_i for pendent}, we can see that
	\begin{align*}
	&\kappa(G)\\
	=&\sum\limits_{i=1}^{k}\frac{m_{G_i}}{m_G}\kappa(G_i)+\sum\limits_{i=1}^{k-1}\left(\frac{m_{G_{i+1}}\hspace*{-0.1cm}+1}{m_G}\mathbf{d}_{G_1}^TR_{G_1}\be_{x_i}\hspace*{-0.1cm}+\frac{m_G-m_{G_{i+1}}}{m_G}\mathbf{d}_{G_{i+1}}^TR_{G_{i+1}}\be_{y_{i}}\right)\\
	&+\hspace*{-0.2cm}\sum\limits_{1\leq i<j\leq k-1}\hspace*{-0.2cm}\frac{2(m_{G_{i+1}}\hspace*{-0.1cm}+1)(m_{G_{j+1}}\hspace*{-0.1cm}+1)}{m_G}r_{x_i,x_j}^{G_1}+\sum\limits_{i=1}^{k-1}\frac{(2m_{G_{i+1}}+1)(2m_G-2m_{G_{i+1}}+1)}{2m_G}.
	\end{align*}
\end{example}

We now recast the formula in \Cref{Theorem:MainResult} in order to provide an intuition for understanding how $\kappa(G)$ is affected by $G_1,\dots,G_k$ and the placement of bridges in $\cB_G$, by expressing the formula for $\kappa(G)$ in terms of the values of $\kappa(G_i)$, accessibility indices, mean first passage times, and the numbers of edges in certain subgraphs of $G$.

\begin{theorem}\label{main 2}
	Let $\cT$ be a tree on $k$ vertices, and $G_1,\dots,G_k$ be connected graphs. Let $G$ be a chain of connected graphs $G_1,\dots,G_k$ with respect to $\cT$. Then, 
	\begin{align}\label{3;formula K(G)}
	\begin{split}
	\kappa(G)  = &\sum_{i=1}^k \kappa(G_i) + \sum_{x\sim y \in \cB_G}\left( \frac{\overline{W}_x}{m_G}\alpha_{G_x}(x) + \frac{\overline{W}_y}{m_G}\alpha_{G_y}(y)\right) \\
	&+\sum\limits_{i=1}^{k}\sum\limits_{\small(z_1,z_2) \in \widetilde{V}_i\times\widetilde{V}_i}\hspace*{-0.2cm}\frac{\overline{W}_{z_1}\overline{W}_{z_2}}{2m_G^2}(m_{z_1, z_2}^{(G_i)}+m_{z_2, z_1}^{(G_i)}) + \hspace*{-0.3cm}\sum\limits_{x\sim y \in \cB_G}\hspace*{-0.3cm}\frac{(2\widebar{W}_x-1)(2\widebar{W}_y-1)}{2m_G},
	\end{split}
	\end{align}
	where $m_{z_1,z_2}^{(G_i)}$ is the mean first passage time from $z_1$ to $z_2$ for a random walk on $G_i$.
\end{theorem}
\begin{proof}
	We recall that for a connected graph $G$, $r_{i,j}^G=\frac{1}{2m_G}(m_{i,j}+m_{j,i})$ where $m_{i,j}$ is the mean first passage time from $i$ to $j$. Then, the conclusion follows from \Cref{Theorem:MainResult} with Lemma \ref{lemma:acce moment kemeny}.
\end{proof}

\subsection{Optimization of Kemeny's constant for chains of connected graphs with respect to trees}

Here we consider how we can maximize/minimize Kemeny's constant for chains of connected graphs with respect to trees, as in the spirit of Section \ref{Subsec:optimization1}. Let $k\geq 2$. Suppose that $\cT$ is a tree on $k$ vertices and $G_1,\dots,G_k$ are $k$ connected graphs. Consider a chain $G$ of $G_1,\dots, G_k$ with respect to $\cT$. Then, for each $x\sim y\in\cB_G$, $x\in V(G_i)$ and $y\in V(G_j)$ for some $i$ and $j$ with $i\neq j$. Considering how $\overline{W}_x$ and $\overline{W}_y$ are defined, we can find that the quantities $\overline{W}_x$ and $\overline{W}_y$ do not depend on the choices of $x$ in $G_i$ and $y$ in $G_j$; further, they only rely on the choice of $\cT$ and $G_1,\dots,G_k$. Therefore, when a tree $\cT$ and connected graphs $G_1,\dots,G_k$ are given, maximizing/minimizing $\kappa(G)$ is equivalent to maximizing/minimizing the second and third summands of the right side of \eqref{3;formula K(G)}.

When it comes to minimization problem, we have fewer constraints. So, we pose and address two minimization problems. 
\begin{question}\label{min-max question1} Given a tree $\cT$ on $k\geq 2$ vertices, and connected graphs $G_1,\dots,G_k$, what is the minimum value of $\kappa(G)$, where $G$ is a chain of $G_1,\dots,G_k$ with respect to $\cT$? 
\end{question}
Choose a chain $G'$ of $G_1,\dots,G_k$ with respect to $\cT$ such that for $i=1,\dots,k$, each $\widetilde{V}_i^{G'}$ (the set of vertices in $G_i$ incident with a bridge in $\cB_{G'}$) is a singleton, and the element $z'$ in $\widetilde{V}_i^{G'}$ satisfies $\alpha_{G_i}(z')\leq \alpha_{G_i}(z)$ for $z\in V(G_i)$. That is, we suppose that only one vertex from each $G_i$ is incident with any bridge in $\cB_{G'}$, and that vertex is one with minimal accessibility index in $G_i$. It follows that $\kappa(G')\leq \kappa(G)$. So, we may ignore the third expression on the right side in the formula \eqref{3;formula K(G)}. Therefore, the minimization problem for Kemeny's constant for all possible chains of $G_1,\dots,G_k$ with respect to $\cT$ is equivalent to a problem of minimizing the accessibility index of a vertex in each of $G_1,\dots,G_k$. 

In \Cref{min-max question1}, in the case that $G_1,\dots,G_k$ are all trees, we are able to produce an answer by using \Cref{centroid minimum}.

\begin{proposition}
	Let $\cT$ be a tree on $k$ vertices, and $G_1,\dots,G_k$ be trees. Then, the minimum of $\kappa(G)$ over all chains $G$ of $G_1,\dots, G_k$ with respect to $\cT$ is attained if and only if for each bridge $x\sim y\in\cB_G$, $x$ and $y$ are centroids of $G_x$ and $G_y$, respectively, and for each $i$, $\widetilde{V}_i^G$ is a singleton.
\end{proposition}

\begin{question}\label{min-max question2} Given a connected graph $H$, what is the minimum value of $\kappa(G)$ for all chains $G$ of $k$ copies of $H$ with respect to a tree $\cT$ on $k$ vertices?
\end{question}
As seen when addressing \Cref{min-max question1}, we may annihilate the third expression of the right side in \eqref{3;formula K(G)} by assuming that exactly one vertex from each $G_i$ is incident with bridges in $\cB_G$. Note that for each $x\sim y\in \cB_G$, $\widebar{W}_x+\widebar{W}_y=m_G+1$. The minimum of the second expression can be obtained by determining a vertex $v$ whose accessibility index in $H$ is minimum, and letting each bridge in $\cB_G$ join the copies of $v$ in each copy of $H$. Now we only need to consider the last summand of the right side in \eqref{3;formula K(G)}. Considering the fact that for $x\sim y\in \cB_G$, $\widebar{W}_x=W_y+1$, $\widebar{W}_y=W_x+1$, and $W_x+W_y=m_G-1$, minimizing Kemeny's constant is equivalent to minimizing the following:
\begin{align}\label{objective function}
\sum\limits_{x\sim y \in \cB_G}W_xW_y.
\end{align}
We can further simplify the question. Note that $W_x$ (resp. $W_y$) can be written in terms of $m_H$ and the number of bridges in $\cB_G$ that belong to the component with $x$ (resp. $y$) in $G\backslash x\sim y$. In the context of minimizing \eqref{objective function}, regarding $m_H$ as $1$, $W_x$ (resp. $W_y$) may be viewed as the sum of the number of vertices and the number of edges in the subtree $\cT_x$ (resp. $\cT_y$) with $x$ (resp. $y$) in $\cT\backslash x\sim y$. Therefore, it follows from the handshaking lemma that we only need to find the minimum of
\begin{align*}
C(\cT):=\sum_{x\sim y \in E(\cT)}|V(\cT_x)| |V(\cT_y)|
\end{align*}
over all trees $\cT$ on $k$ vertices. 

We introduce further definitions and notation to address \Cref{min-max question2}. For a connected graph $G$ with a vertex $v$, if $G-v$ has $r$ connected components $H_1,\dots,H_r$ for some $r\geq 2$, then the subgraph induced by $V(H_i)$ for $1\leq i\leq r$ is called a \textit{branch} of $G$ at $v$. Let $\cT$ be a tree, and $v$ be a vertex of $\cT$. For $w\in V(\cT)\backslash\{v\}$, we use $c_v(w)$ to denote the number of vertices of the subtree obtained from $\cT$ by removing the branch of $\cT$ at $w$ that contains $v$.

\begin{remark}\label{Remark: c_v(w)}
	In this remark, we discuss which tree minimizes the quantity $\sum_{x\in V(\cT)\backslash\{v\}} c_v(x)$, provided $v$ is a pendent vertex of a tree $\cT$ on $n$ vertices. Suppose that $B$ is the branch of $\cT$ at $w$ that contains $v$. Then, $c_v(w)=n-|V(B)|$. Clearly, for $w\in V(\cT)\backslash\{v\}$, $1\leq c_v(w)\leq n-1$. Since $v$ is a pendent vertex, $c_v(w)=n-1$ if and only if $w$ is adjacent to $v$; so, we have exactly one vertex $w$ with $c_v(w)=n-1$. If $w$ is a pendent vertex, then $c_v(w)=1$. It follows that when $v$ is a pendent vertex, $\sum_{x\in V(\cT)\backslash\{v\}} c_v(x)$ is minimized if and only if $\cT$ is a star.
\end{remark}

In order to answer \Cref{min-max question2}, we consider a property of $C(\cT)$ where $\cT$ is a tree on $n$ vertices. Let $v$ be a pendent vertex of $\cT$, and let $\cT'=\cT-v$. For each $x\sim y \in E(\cT)$, either $\cT_x$ or $\cT_y$ contains $v$, so assuming that in the following, $\cT_y$ does so, we obtain
\begin{align}\nonumber
C(\cT)=&\sum_{x\sim y \in E(\cT)}|V(\cT_x)| |V(\cT_y)|\\\nonumber
=&\sum_{x\sim y \in E(\cT')}|V(\cT'_x)| (|V(\cT'_y)|+1) + (n-1) \\\label{expression for induction}
=& \, C(\cT') + \sum_{w\in V(\cT)\backslash\{v\}}c_v(w) +(n-1).
\end{align}
(Note that in \eqref{expression for induction}, $c_v(w)$ is the number of vertices of the subtree obtained from $\cT$ (not $\cT'$) by removing the branch of $\cT$ at $w$ that contains $v$.)

\begin{lemma}
	Let $n\geq 4$. The minimum of $C(\cT)$ for trees $\cT$ on $n$ vertices is attained if and only if $\cT$ is a star.
\end{lemma}
\begin{proof}
	We shall use induction on $n$ for this proof. It can be seen from computation that the statement holds for $n=4$. Let $n\geq 5$. Suppose that $\cT$ is a tree on $n$ vertices. We choose a pendent vertex $v$ in $\cT$. Consider $\cT'=\cT-v$. We can see from the inductive hypothesis and \eqref{expression for induction} that it suffices to show that the minimum of $\sum_{w\in V(\cT)\backslash\{v\}}c_v(w)$ is attained if and only if $\cT$ is a star with a pendent vertex $v$. By \Cref{Remark: c_v(w)}, the conclusion follows.
\end{proof}

The following states the answer to \Cref{min-max question2}.

\begin{proposition}
	Let $H$ be a connected graph. Fix a vertex $v\in V(H)$ such that $\alpha_{H}(v)\leq \alpha_{H}(w)$ for all $w\in V(H)$. Then, the minimum of $\kappa(G)$ over all chains $G$ of $k$ copies of $H$ with respect to any tree $\cT$ on $k$ vertices is attained when $\cT$ is a star and for each bridge $x\sim y\in\cB_G$, $x=v$ and $y=v$.
\end{proposition}

\section{Concluding remarks}\label{Sec4}
While we are devoted to deriving a formula of Kemeny's constant for a graph with bridges in terms of several quantities inherent to the subgraphs obtained upon the deletion of the bridges, we argue that this result may find many interesting applications in practice. Here, we give some suggestions of how our findings could be used in practice, without fully developing the identified applications. 

If one has several connected graphs $G_1,\dots, G_k$, whose resistance matrices and degree vectors are known, then Theorems \ref{Theorem:MainResult} and \ref{main 2} may be used to decide how we should connect them with several edges so that the resulting graph is a chain $G$ of $G_1,\dots,G_k$ with respect to some tree $\cT$, in order to maximize/minimize Kemeny's constant for the resulting graph. In particular, as discussed in \Cref{min-max question1} and \Cref{min-max question2}, when it comes to minimizing Kemeny's constant, there are fewer constraints than the problem of maximizing Kemeny's constant does.

This could be useful when $G_1,\dots,G_k$ correspond to some networks (e.g., microgrids in power systems in emerging countries, or transportation networks, such as networks of air flights), and one has to decide which nodes of the networks should be connected in order to maximize the connectivity of the new whole network, which corresponds to $G$. As Kemeny's constant is known to provide a useful indication on the connectivity of a transportation network (see for instance \cite{Crisostomi:Google}), then in the second example the objective could be to find which two airports of two different flight networks should be connected in order to minimize the Kemeny's constant of the overall transportation network (and thus, maximize its connectivity).

Besides, there is some interesting application regarding maximizing Kemeny's constant. That is, we shall connect $G_1,\dots,G_k$ with several edges so that the resulting graph is a chain $G$ of $G_1,\dots,G_k$ with respect to some tree, while keeping them as least connected as possible, and thus maximizing $\kappa(G)$. Related examples have recently emerged in the pandemic scenario, when one may be interested in connecting social/community networks in the least connected possible way. For instance, if $G_1,\dots,G_k$ correspond to different and non-connected departments inside companies, the corresponding heads of departments may want to plan meetings to mitigate the risk of spreading the virus between departments so that each head has at least one meeting and the number of events is $k-1$; hence, one can arrange schedules in such a way to maximize $\kappa(G)$.

Finally, Theorem \ref{Theorem:MainResult} may be used to compute Kemeny's constant of a large graph with bridges in a parallel fashion. We note Remark \ref{Remark:complexity}. If one understand where those bridges are in the graph, and if the resulting graph after deleting the bridges contains connected components `similar' in size, then the theorem can be used to efficiently compute Kemeny's constant for the large graph, by understanding the resistance matrices and degree vectors of the connected components.

\noindent{\em Acknowledgements}: E.C. and S.K. were supported by the Research Project PRIN 2017 “Advanced Network Control of Future Smart Grids” funded by the Italian Ministry of University and Research (2020–2023).\\
J.B. was supported by NSERC Discovery Grant RGPIN-2021-03775.


\end{document}